\documentclass[11pt,a4paper]{article}
\usepackage[T1]{fontenc}
\usepackage[utf8]{inputenc}
\usepackage{lmodern,amsmath,amssymb,amsthm,mathtools}
\usepackage[margin=27mm]{geometry}
\usepackage{enumitem,microtype}
\usepackage[hidelinks]{hyperref}
\newtheorem{theorem}{Theorem}[section]
\newtheorem{proposition}[theorem]{Proposition}
\newtheorem{lemma}[theorem]{Lemma}
\newtheorem{corollary}[theorem]{Corollary}
\theoremstyle{remark}
\newtheorem{remark}[theorem]{Remark}
\newtheorem{example}[theorem]{Example}
\newcommand{\B}{\mathcal B}
\newcommand{\C}{\mathbb C}
\newcommand{\rad}{\operatorname{r}}
\newcommand{\codim}{\operatorname{codim}}
\setlist[enumerate,1]{label=\textup{(\roman*)},leftmargin=9mm}
\title{Independent Backward Preimages and\\Algebraic Radial Decompositions of Operators}
\author{Zakariae Aznay\textsuperscript{1}\quad
Fadil Chabbabi\textsuperscript{1}\quad
Hassan Zariouh\textsuperscript{2}\\[8pt]
\normalsize\textsuperscript{1}Laboratory (A.F.N.L.A.), Department of Mathematics,\\[-1pt]
\normalsize Faculty of Sciences, Abdelmalek Essa\^adi University,\\[-1pt]
\normalsize Tetouan, Morocco\\[2pt]
\normalsize E-mail (Z. Aznay): \texttt{aznay.zakariae@gmail.com}\\[-1pt]
\normalsize E-mail (F. Chabbabi): \texttt{f.chabbabi@uae.ac.ma}\\[7pt]
\normalsize\textsuperscript{2}Laboratory LCSI2A, Faculty of Sciences,\\[-1pt]
\normalsize Mohammed First University, Oujda, Morocco\\[2pt]
\normalsize E-mail (H. Zariouh): \texttt{h.zariouh@yahoo.fr}}
\date{}
\begin{document}
\maketitle

\begin{abstract}
Let $T\in\mathcal{B}(X)$, where $X$ is a complex Banach space.
We introduce the spaces $\mathcal{I}_s(T)$, defined by independently
chosen preimages under the powers of $T$ with exponential growth
bounds, and study algebraic radial decompositions involving
$\mathcal{H}_r(T)$, $\mathcal{K}_s(T)$, and $\mathcal{I}_s(T)$.
We prove that every algebraic decomposition
$X=\mathcal{H}_r(T)\oplus\mathcal{I}_s(T)$ is topological and that
$\mathcal{I}_s(T)=\mathcal{K}_s(T)$.
This removes the closedness assumption in the radial decomposition
theorem of Koliha and Poon.
We also prove that, for $0\leq a\leq b$, the closedness of
$\mathcal{I}_b(T)$ and
$\mathcal{H}_a(T)+\mathcal{I}_b(T)$ implies
$\mathcal{I}_b(T)=\mathcal{K}_b(T)$.
We establish finite-defect and circular spectral criteria and
characterize the Cantor--Bendixson derivatives of the spectrum,
Browder spectrum, and Drazin spectrum in terms of radial
decompositions and quotient restrictions.
As applications, we obtain new characterizations of generalized
Drazin, generalized Drazin--Riesz, generalized
Drazin--meromorphic, and $g_z$-invertible operators
in terms of radial subspaces.

\end{abstract}



\section{Introduction}
Let $X$ be a complex Banach space and let $\B(X)$ denote the algebra of all bounded linear operators on $X$. For $T\in\B(X)$, let $\mathcal R(T)$ and $\mathcal N(T)$ denote its range and kernel, respectively. For $r>0$, put
\begin{equation*}
\mathcal{H}_r(T)=\{x\in X:\limsup_{k\to\infty}\|T^kx\|^{1/k}<r\}
\end{equation*}
and
\begin{equation*}
\mathcal{K}_r(T)=\left\{x\in X:\begin{array}{l}
\exists(x_k)_{k\geq0}\subset X,\quad x_0=x,\quad Tx_{k+1}=x_k,\\
\limsup_{k\to\infty}\|x_k\|^{1/k}<r^{-1}
\end{array}\right\}.
\end{equation*}
For $r=0$, the quasinilpotent part and the analytic core are given by
\begin{align}
\mathcal{H}_0(T)&=\{x\in X:\limsup_{k\to\infty}\|T^kx\|^{1/k}=0\}
       =\bigcap_{r>0}\mathcal{H}_r(T),\notag\\
\mathcal{K}(T)&=\left\{x\in X:\begin{array}{l}
\exists(x_k)_{k\geq0}\subset X,\ x_0=x,\ Tx_{k+1}=x_k,\\
\limsup_{k\to\infty}\|x_k\|^{1/k}<\infty
\end{array}\right\}
       =\bigcup_{r>0}\mathcal{K}_r(T).\label{zero}
\end{align}
We write $\mathcal K_0(T)=\mathcal K(T)$.

Koliha and Poon studied the radial spaces $\mathcal H_r(T)$ and $\mathcal K_r(T)$ in \cite[Section 5]{KP}. At radius zero these spaces are the quasinilpotent part and the analytic core \cite{Mbekhta}. The algebraic equality $X=\mathcal H_0(T)\oplus\mathcal K(T)$ already implies generalized Drazin invertibility \cite[Theorem 3.1]{AC}. We consider arbitrary nonnegative radii and compare coherent backward chains with independently chosen preimages.

We introduce, for $r>0$, the subspaces
\[
\mathcal I_r(T)=\left\{x\in X:\begin{array}{l}
\exists(x_n)_{n\geq0}\subset X,\quad T^nx_n=x\quad(n\geq0),\\
\displaystyle\limsup_{n\to\infty}\|x_n\|^{1/n}<r^{-1}
\end{array}\right\},
\]
and put
\[
\mathcal I(T)=\left\{x\in X:\begin{array}{l}
\exists(x_n)_{n\geq0}\subset X,\quad T^nx_n=x\quad(n\geq0),\\
\displaystyle\limsup_{n\to\infty}\|x_n\|^{1/n}<\infty
\end{array}\right\},\qquad \mathcal I_0(T)=\mathcal I(T).
\]
For $n\geq1$, we equip $\mathcal R(T^n)$ with the range norm
\[
\|x\|_{T^n}=\inf\{\|y\|:T^ny=x\},\qquad x\in\mathcal R(T^n),
\]
and put $\|x\|_{T^n}=+\infty$ when $x\notin\mathcal R(T^n)$.
This is the quotient norm induced by $X/\mathcal N(T^n)$, so
$\mathcal R(T^n)$ is complete in this norm even when it is not closed in $X$.
Equivalently,
\begin{align*}
\mathcal I_r(T)&=\left\{x\in X:\limsup_{n\to\infty}\|x\|_{T^n}^{1/n}<r^{-1}\right\}\quad(r>0),\\
\mathcal I(T)&=\left\{x\in X:\limsup_{n\to\infty}\|x\|_{T^n}^{1/n}<\infty\right\}.
\end{align*}

A backward chain supplies independent preimages at every order, so $\mathcal K_b(T)\subseteq\mathcal I_b(T)$. The reverse inclusion can fail even when $\mathcal I_b(T)$ is one-dimensional. Theorem~\ref{closedindcores} gives conditions under which every vector in $\mathcal I_b(T)$ admits a backward chain with the prescribed growth. Theorem~\ref{finitecircle} identifies the defect of a radial sum at a fixed radius, and Theorem~\ref{radialtransfinite} treats finite defects at two possibly different radii.

We denote by $\sigma(T)$ and $\rho(T)$ the spectrum and resolvent set of $T$, respectively, and by $\rad(T)$ its spectral radius. The Browder and Drazin spectra are
\[
\begin{aligned}
\sigma_b(T)&=\sigma(T)\setminus\{\lambda:\lambda\text{ is a pole of finite rank of the resolvent}\},\\
\sigma_d(T)&=\sigma(T)\setminus\{\lambda:\lambda\text{ is a pole of the resolvent}\}.
\end{aligned}
\]
For $E\subseteq\C$, let $\operatorname{acc}E$ and $\operatorname{iso}E$ denote its accumulation and isolated points, respectively. Put $\mathbb T_r=\{\lambda\in\C:|\lambda|=r\}$. We write $X=M\oplus N$ for an algebraic direct sum; no closedness is assumed. If $M$ and $N$ are closed, the direct sum is topological. For a $T$-invariant subspace $M$, we write $T_M:M\to M$ for the restriction of $T$ to $M$, endowed with the norm inherited from $X$. Zero summands are allowed; the spectrum of an operator on $\{0\}$ is empty and its spectral radius is taken to be zero. We use $I$ for the identity on the space in question and $\codim L=\dim(X/L)$ for algebraic codimension. A reducing decomposition consists of two closed invariant subspaces whose direct sum is $X$. An operator is generalized Drazin invertible if it has a reducing decomposition into an invertible and a quasinilpotent operator \cite[Theorem 7.1]{Koliha}.

\section{Radial subspaces}
For $r>0$, the spaces $\mathcal{H}_r(T)$ and $\mathcal{K}_r(T)$ are hyperinvariant. Shifting forward or backward iterates gives
\begin{equation}\label{inv}
T^{-1}(\mathcal{H}_r(T))=\mathcal{H}_r(T),\qquad T\mathcal{K}_r(T)=\mathcal{K}_r(T),
\end{equation}
where $T^{-1}$ denotes preimage. Moreover,
\begin{equation*}
0<s<r\quad\Longrightarrow\quad \mathcal{H}_s(T)\subset \mathcal{H}_r(T),\quad \mathcal{K}_r(T)\subset \mathcal{K}_s(T).
\end{equation*}
The identities in \eqref{inv} also hold at radius zero.

Here $c_0(X)$ denotes the Banach space of norm-null sequences in $X$, equipped with the supremum norm $\|\cdot\|_\infty$. For $a,b>0$, define
\begin{align*}
V_a&=\{x\in X:a^{-k}T^kx\longrightarrow0\},&
\|x\|_a&=\sup_{k\geq0}a^{-k}\|T^kx\|,\\
Y_b&=\{y=(y_k)\in c_0(X):bTy_{k+1}=y_k\ (k\geq0)\},&
E_b(y)&=y_0.
\end{align*}
\begin{proposition}\label{models}
The spaces $(V_a,\|\cdot\|_a)$ and $Y_b$ are Banach spaces. The inclusion $V_a\to X$ and the map $E_b:Y_b\to X$ are bounded. For every $r>0$,
\begin{equation*}
\mathcal{H}_r(T)=\bigcup_{\substack{a\in\mathbb Q\\0<a<r}}V_a,
\qquad
\mathcal{K}_r(T)=\bigcup_{\substack{b\in\mathbb Q\\0<b<r^{-1}}}E_b(Y_b).
\end{equation*}
\end{proposition}
\begin{proof}
The maps
\[
x\longmapsto(a^{-k}T^kx)_{k\geq0},\qquad
 y\longmapsto(bTy_{k+1}-y_k)_{k\geq0}
\]
show that $V_a$ is isometric to a closed subspace of $c_0(X)$ and that $Y_b$ is closed in $c_0(X)$. The bounds $\|x\|\leq\|x\|_a$ and $\|E_b y\|\leq\|y\|_\infty$ give continuity. If $x\in \mathcal{H}_r(T)$, choose a rational $a$ with
$\limsup_{k\to\infty}\|T^kx\|^{1/k}<a<r$.
Then $a^{-k}T^kx\to0$, so $x\in V_a$. Conversely, if $x\in V_a$ for some rational $0<a<r$, then
$\limsup_{k\to\infty}\|T^kx\|^{1/k}\leq a<r$.

If $x\in \mathcal{K}_r(T)$, let $(x_k)$ be a determining chain and choose a rational $b$ such that
\[
\limsup_{k\to\infty}\|x_k\|^{1/k}<b<r^{-1}.
\]
Putting $y_k=b^{-k}x_k$ gives $y\in Y_b$ and $E_b(y)=x$. Conversely, if $x=E_b(y)$ for some $y\in Y_b$ and rational $0<b<r^{-1}$, then $x_k=b^ky_k$ satisfies $x_0=x$, $Tx_{k+1}=x_k$, and
\[
\limsup_{k\to\infty}\|x_k\|^{1/k}\leq b<r^{-1}.
\]
\end{proof}

\smallskip
\noindent By \cite[Proposition 2.3]{AC}, $\mathcal{H}_0(T)$ is a Fr\'echet space under the norms
\[
p_m(x)=\sup_{k\geq0}m^k\|T^kx\|\qquad(m\geq1),
\]
and its inclusion into $X$ is continuous. Also, \eqref{zero} gives
\[
\mathcal{K}_0(T)=\bigcup_{b\in\mathbb Q,\ b>0}E_b(Y_b).
\]
We shall use the nonmeagre-range form of the open mapping theorem \cite[Theorem 2.11]{Rudin}.

\begin{proposition}\label{dichotomy}
Let $T\in\B(X)$ and $r,s\geq0$, and put $E=\mathcal{H}_r(T)\cap \mathcal{K}_s(T)$, where $\mathcal{K}_0(T)=\mathcal{K}(T)$. If $0<\dim E<\infty$, then $T_E$ is invertible and
\[
\sigma(T_E)\subset\{\lambda\in\C:s<|\lambda|<r\}.
\]
Consequently, for $0\leq r\leq s$, the intersection $\mathcal{H}_r(T)\cap \mathcal{K}_s(T)$ is either zero or infinite-dimensional.
\end{proposition}

\begin{proof}
Let $0<\dim E<\infty$. If $x\in E$ and $(x_n)$ is a determining chain for $x\in \mathcal{K}_s(T)$, then $x_n\in \mathcal{K}_s(T)$ and, by \eqref{inv},
\[
T^nx_n=x\in \mathcal{H}_r(T)\quad\Longrightarrow\quad x_n\in \mathcal{H}_r(T).
\]
Thus $TE=E$ and $T_E$ is invertible. For any eigenvalue $\lambda$ of $T_E$, take $0\ne x\in E$ with $Tx=\lambda x$. Then $x_n=\lambda^{-n}x$ and
\[
|\lambda|=\lim_{n\to\infty}\|T^nx\|^{1/n}<r,
\qquad
|\lambda|^{-1}=\lim_{n\to\infty}\|x_n\|^{1/n}<s^{-1}
\]
when $r,s>0$. If $r=0$, then $x\in \mathcal{H}_0(T)$ gives $|\lambda|=0$, contradicting invertibility. If $s=0$, invertibility gives $|\lambda|>0$.
\end{proof}

\section{Independent backward preimages}\label{sec:independent}
\begin{proposition}\label{indbasic}
Let $r>0$ and $m\geq1$. The spaces $\mathcal I_r(T)$ and
$\mathcal I(T)$ are hyperinvariant linear subspaces, and
\[
\mathcal K_r(T)\subseteq\mathcal I_r(T)\subseteq\mathcal R^\infty(T),
\qquad
\mathcal K(T)\subseteq\mathcal I(T)\subseteq\mathcal R^\infty(T),
\]
where $\mathcal R^\infty(T)=\bigcap_{n\geq1}\mathcal R(T^n)$. Moreover,
\[
0<s<r\ \Longrightarrow\ \mathcal I_r(T)\subseteq\mathcal I_s(T),
\qquad
\mathcal I(T)=\bigcup_{r>0}\mathcal I_r(T),
\]
and
\begin{equation}\label{indpowers}
\mathcal I_r(T)=\mathcal I_{r^m}(T^m),
\qquad \mathcal I(T)=\mathcal I(T^m).
\end{equation}
If $T$ is injective, then $\mathcal I_r(T)=\mathcal K_r(T)$
and $\mathcal I(T)=\mathcal K(T)$. For a topological invariant
direct sum $T=A\oplus B$,
\[
\mathcal I_r(T)=\mathcal I_r(A)\oplus\mathcal I_r(B),
\qquad\mathcal I(T)=\mathcal I(A)\oplus\mathcal I(B).
\]
\end{proposition}
\begin{proof}
The inclusions and hyperinvariance follow by applying an operator
commuting with $T$ to the determining preimages. If $T^nx_n=x$,
then $T^{mn}x_{mn}=x$; this proves one inclusion in
\eqref{indpowers}. Conversely, if $T^{mn}y_n=x$, set
$z_{mn}=y_n$ and $z_{mn+j}=T^{m-j}y_{n+1}$ for $1\leq j<m$.
Then $T^kz_k=x$ for all $k$. Since
$\limsup_{n\to\infty}\|y_n\|^{1/n}<r^{-m}$, choose
$0<q<r^{-m}$ and $C>0$ such that $\|y_n\|\leq Cq^n$ for all
$n\geq0$. As $k=mn+j$ with $0\leq j<m$, there exists $C'>0$
such that $\|z_k\|\leq C'q^{k/m}$ for every $k\geq0$. Hence
\[
\limsup_{k\to\infty}\|z_k\|^{1/k}\leq q^{1/m}<r^{-1}.
\]
The same construction gives $\mathcal I(T)=\mathcal I(T^m)$,
since it preserves finite exponential growth.
If $T$ is injective, $T^n(Tx_{n+1}-x_n)=0$ implies
$Tx_{n+1}=x_n$. For a topological reducing decomposition, applying the bounded
projections to preimages gives the corresponding preimages in each
summand; taking their direct sum gives the converse.
\end{proof}

\begin{lemma}\label{indfinitecodim}
Let $T\in\B(X)$ and let $L$ be a closed $T$-invariant subspace
of finite codimension. Then, for every $b\geq0$,
\[
\mathcal I_b(T_L)=\mathcal I_b(T)\cap L.
\]
\end{lemma}
\begin{proof}
We need only prove that $\mathcal I_b(T)\cap L\subseteq\mathcal I_b(T_L)$.
Let $\widetilde T$ be induced by $T$ on $X/L$. Since $X/L$ is
finite-dimensional, choose $p\geq0$ such that
\[
\mathcal N(\widetilde T^{n+p})=\mathcal N(\widetilde T^p)
\qquad(n\geq0).
\]
Let $x\in\mathcal I_b(T)\cap L$, and choose $x_n\in X$ with
$T^nx_n=x$ and the prescribed exponential growth. Since
$T^{n+p}x_{n+p}=x\in L$, the equality of the kernels gives
$T^px_{n+p}\in L$. Put $y_n=T^px_{n+p}$. Then
\[
T_L^ny_n=x,\qquad
\|y_n\|\leq\|T^p\|\,\|x_{n+p}\|,
\]
and hence
\[
\limsup_{n\to\infty}\|y_n\|^{1/n}
\leq\limsup_{n\to\infty}\|x_n\|^{1/n}.
\]
Thus $x\in\mathcal I_b(T_L)$.
\end{proof}

\begin{remark}\label{indradremark}
If $0\ne x\in\mathcal R^\infty(T)$, then
$\|x\|\leq\|T^n\|\,\|x\|_{T^n}$ for every $n\geq1$. Hence
\[
\liminf_{n\to\infty}\|x\|_{T^n}^{1/n}
\geq\rad(T)^{-1},
\]
with $1/0=+\infty$. Consequently,
$\mathcal I_r(T)=\{0\}$ whenever $\rad(T)\leq r$, for $r\geq0$.
\end{remark}

\smallskip
\noindent For $b>0$, define
\[
F_b=\{z=(z_n)\in c_0(X):T^n(b^nz_n)=z_0\ (n\geq1)\},
\qquad J_b(z)=z_0.
\]
\begin{proposition}\label{indmodels}
Each $F_b$ is a Banach space and $J_b:F_b\to X$ is bounded.
For every $r>0$,
\[
\mathcal I_r(T)=\bigcup_{\substack{b\in\mathbb Q\\0<b<r^{-1}}}
\mathcal R(J_b),\qquad
\mathcal I(T)=\bigcup_{\substack{b\in\mathbb Q\\b>0}}
\mathcal R(J_b).
\]
Furthermore,
\begin{align*}
\mathcal I_r(T)\text{ is nonmeagre in }X
&\iff\mathcal I_r(T)=X\iff\mathcal K_r(T)=X,\\
\mathcal I(T)\text{ is nonmeagre in }X
&\iff\mathcal I(T)=X\iff T\text{ is surjective}.
\end{align*}
\end{proposition}
\begin{proof}
Clearly, $F_b$ is a closed subspace of $c_0(X)$ and $J_b$
is bounded. The range representations follow by taking
$z_n=b^{-n}x_n$.
Suppose $\mathcal I_r(T)$ is nonmeagre. By the countable
range representation above, some $\mathcal R(J_b)$ with
$b<r^{-1}$ is nonmeagre. By \cite[Theorem 2.11]{Rudin},
$J_b$ is onto and open. Thus
there exists $C>0$ such that
\[
\forall x\in X,\ \forall n\geq1,\quad
\exists y\in X:\quad T^ny=x,\quad \|y\|\leq Cb^n\|x\|.
\]
Choose $m\geq1$ so that $(Cb^m)^{1/m}<r^{-1}$.
Repeatedly choosing preimages under $T^m$ gives
$u_0=x$, $T^mu_{k+1}=u_k$ and
$\|u_k\|\leq(Cb^m)^k\|x\|$. Defining
\[
x_{km}=u_k,\qquad x_{km+j}=T^{m-j}u_{k+1}
\quad(1\leq j<m)
\]
yields $Tx_{n+1}=x_n$ and
$\limsup_{n\to\infty}\|x_n\|^{1/n}
\leq(Cb^m)^{1/m}<r^{-1}$. Thus $\mathcal K_r(T)=X$.
The converse is immediate. Since $\mathcal I(T)\subseteq\mathcal R(T)$,
nonmeagreness of $\mathcal I(T)$ implies that $T$ is
surjective by \cite[Theorem 2.11]{Rudin}. Conversely,
$\mathcal K(T)=X$ for surjective $T$ by
\cite[Lemma 3(5)]{GMO}.
\end{proof}

\smallskip
\noindent The inclusions in Proposition \ref{indbasic} may be strict, even on a
Hilbert space. In particular, $T\mathcal I_r(T)$ need not equal
$\mathcal I_r(T)$.
\begin{example}\label{indcounterexample}
Let $\mathcal H$ have the orthonormal basis
$\{e_0\}\cup\{e_{m,k}:m\geq1,\ 1\leq k\leq m\}$,
and define $T\in\B(\mathcal H)$ by
\[
Te_0=0,\qquad Te_{m,1}=4\,2^{-m}e_0,\qquad
Te_{m,k}=4e_{m,k-1}\quad(2\leq k\leq m).
\]
The operator is bounded since $(4\,2^{-m})_{m\ge1}\in\ell^2$ and all other nonzero weights are $4$. We have $T^ne_{n,n}=2^ne_0$, so $e_0\in\mathcal I_1(T)$.
For $n>m-k$, the $(m,k)$-coordinate of every vector in $\mathcal R(T^n)$ vanishes. Hence $\mathcal R^\infty(T)=\C e_0$ and
$T e_0=0$. If $(y_n)_{n\geq0}$ is a backward chain, then
$y_n=T^ky_{n+k}\in\mathcal R^\infty(T)=\C e_0$
for all $n,k\geq0$. Since $Te_0=0$, we have
$y_n=Ty_{n+1}=0$ for every $n\geq0$. Hence
$\mathcal K(T)=\{0\}$. Consequently,
\[
\mathcal I_1(T)=\mathcal I(T)=\C e_0
\subseteq\mathcal H_0(T),\qquad T\mathcal I_1(T)=\{0\}.
\]
\end{example}

\smallskip
\noindent Write $\mathcal C(T)$ for the algebraic core, the largest subspace $M$ satisfying $TM=M$. Recall that $\mathcal K(T)\subseteq\mathcal C(T)\subseteq\mathcal R^\infty(T)$ \cite{Mbekhta}. Neither $\mathcal I(T)$ nor $\mathcal C(T)$ contains the other in general. Indeed, in Example \ref{indcounterexample}, $e_0\in\mathcal I(T)$ but $\mathcal C(T)=\{0\}$. Conversely, for $De_j=j^{-1}e_j$ on $\ell^2$ and $x=\sum_{j\geq1}e^{-j^2}e_j$, the vectors $x_n=\sum_{j\geq1}j^ne^{-j^2}e_j$ form a backward chain, so $x\in\mathcal C(D)$. They are the unique preimages of $x$ under $D^n$, and $\|x_n\|\ge j^ne^{-j^2}$ for each fixed $j$, so $\liminf_n\|x_n\|^{1/n}=\infty$. Thus $x\notin\mathcal I(D)$.

\medskip
\noindent Write $\sigma_{su}(T)=\{\lambda:\mathcal R(T-\lambda I)\ne X\}$. Gonz\'alez, Mbekhta and Oudghiri characterized the isolated
points of this spectrum using $\mathcal H_0(T)+\mathcal K(T)$
\cite[Theorem 5]{GMO}. We characterize $0\notin\operatorname{acc}\sigma_{su}(T)$ using independent preimages and fixed radial bounds.
\begin{theorem}\label{indsurjspectrum}
Let $T\in\B(X)$. Then
\begin{align*}
0\notin\operatorname{acc}\sigma_{su}(T)
&\iff X=\mathcal H_0(T)+\mathcal I(T)\\
&\iff X=\mathcal H_0(T)+\mathcal I_b(T)
\quad\text{for some }b>0.
\end{align*}
These conditions are also equivalent to
$X=\mathcal H_0(T)+\mathcal K_b(T)$ for some $b>0$.
\end{theorem}
\begin{proof}
Suppose first that $0\notin\operatorname{acc}\sigma_{su}(T)$.
If $T$ is surjective, then $\mathcal K(T)=X$ by
\cite[Lemma 3(5)]{GMO}. Otherwise, $0$ is isolated in $\sigma_{su}(T)$,
and \cite[Theorem 5]{GMO} gives
$X=\mathcal H_0(T)+\mathcal K(T)$.
By Proposition \ref{models},
$\mathcal K(T)=\bigcup_{c\in\mathbb Q,\ c>0}E_c(Y_c)$.
Applying the nonmeagre-range open mapping theorem to the maps
\[
\mathcal H_0(T)\times Y_c\longrightarrow X,\qquad
(h,z)\longmapsto h+E_cz,
\]
we obtain $X=\mathcal H_0(T)+E_c(Y_c)$ for some $c>0$.
Since $E_c(Y_c)\subseteq\mathcal K_b(T)$ whenever
$0<b<c^{-1}$, this gives
$X=\mathcal H_0(T)+\mathcal K_b(T)$ for some $b>0$.
The inclusions $\mathcal K_b(T)\subseteq\mathcal I_b(T)
\subseteq\mathcal I(T)$ give the corresponding implications.

Conversely, suppose $X=\mathcal H_0(T)+\mathcal I(T)$.
For every integer $b\geq1$, let $F_b$ and $J_b$ be as in
Proposition \ref{indmodels}, and define
\[
\Phi_b:\mathcal H_0(T)\times F_b\longrightarrow X,\qquad
\Phi_b(h,z)=h+J_bz.
\]
Since $\mathcal I(T)=\bigcup_{b\geq1}\mathcal R(J_b)$,
the ranges of the continuous maps $\Phi_b$ cover $X$.
By the nonmeagre-range open mapping theorem, some $\Phi_b$
is an open surjection.
For each $m\geq1$, there is $C_m>0$ such that every $x\in X$
can be written as
$x=h+z_0$, where $h\in\mathcal H_0(T)$ and $z\in F_b$, with
\begin{equation*}
p_m(h)+\|z\|_\infty\leq C_m\|x\|.
\end{equation*}

Fix $0<|\lambda|<b^{-1}$ and choose $m\geq1$ such that
$m^{-1}<|\lambda|$. For $x=h+z_0$ as above, put
$v_n=b^nz_n$. Then
\[
T^nv_n=z_0,\qquad \|v_n\|\leq b^n\|z\|_\infty.
\]
Since $\|T^jh\|\leq m^{-j}p_m(h)$ and $m^{-1}<|\lambda|$,
the series
\[
y_h=-\sum_{j=0}^\infty\lambda^{-j-1}T^jh
\]
converges in $X$, with
\[
(T-\lambda I)y_h=h,\qquad
\|y_h\|\leq\frac{p_m(h)}{|\lambda|-m^{-1}}.
\]
For $n\geq1$, put
\[
y_{z,n}=\sum_{j=0}^{n-1}\lambda^jT^{n-1-j}v_n.
\]
The identity
\[
(T-\lambda I)\sum_{j=0}^{n-1}\lambda^jT^{n-1-j}
=T^n-\lambda^nI
\]
gives
\[
(T-\lambda I)y_{z,n}=z_0-\lambda^nv_n,\qquad
\|y_{z,n}\|\leq b^n\|z\|_\infty
\sum_{j=0}^{n-1}|\lambda|^j\|T\|^{n-1-j}.
\]
Consequently,
\[
(T-\lambda I)(y_h+y_{z,n})=x-\lambda^nv_n.
\]
Since $|\lambda|b<1$, choose $n\geq1$ such that
$C_m(|\lambda|b)^n\leq\frac12$. The preceding estimates
show that there is $D_\lambda>0$ such that, for every $x\in X$,
one can find $y\in X$ satisfying
\[
\|x-(T-\lambda I)y\|\leq\tfrac12\|x\|,\qquad
\|y\|\leq D_\lambda\|x\|.
\]
Starting with $r_0=x$, construct sequences $(r_k)$ and $(y_k)$
so that
\[
r_{k+1}=r_k-(T-\lambda I)y_k,\qquad
\|r_k\|\leq2^{-k}\|x\|,\quad
\|y_k\|\leq D_\lambda2^{-k}\|x\|.
\]
Then $y=\sum_{k=0}^\infty y_k$ converges in $X$, and
\[
(T-\lambda I)y
=\lim_{k\to\infty}(x-r_{k+1})=x.
\]
Thus $T-\lambda I$ is surjective for $0<|\lambda|<b^{-1}$,
and $0\notin\operatorname{acc}\sigma_{su}(T)$.
\end{proof}

The next two results give stability properties of the radial subspaces under nilpotent and quasinilpotent perturbations.
\begin{theorem}\label{indnilpotent}
For $S,T\in\B(X)$ and $n\geq0$, put
\[
D_n(S,T)=\sum_{k=0}^n(-1)^k\binom nk S^{n-k}T^k.
\]
Suppose that $D_p(S,T)=D_q(T,S)=0$ for some $p,q\geq1$.
Then, for every $x\in X$,
\[
\limsup_{n\to\infty}\|x\|_{S^n}^{1/n}
=\limsup_{n\to\infty}\|x\|_{T^n}^{1/n}.
\]
Consequently, $\mathcal I_r(S)=\mathcal I_r(T)$ for all $r\geq0$.
In particular, these conclusions hold for $S=T+N$ whenever
$N$ is nilpotent and $N\in\operatorname{comm}_w(T)$, where
$\operatorname{comm}_w(T)$ is defined by
$TNT=T^2N=NT^2$ and $TN^2=NTN=N^2T$; see \cite{GCP}.
\end{theorem}
\begin{proof}
We first show that, for every $N\geq0$,
\[
T^N=\sum_{k=0}^{N}(-1)^k\binom Nk S^{N-k}D_k(S,T).
\]
Indeed, using $\binom Nk\binom kj=\binom Nj\binom{N-j}{k-j}$,
\[
\begin{aligned}
\sum_{k=0}^{N}(-1)^k\binom Nk S^{N-k}D_k(S,T)
&=\sum_{j=0}^{N}\binom Nj S^{N-j}T^j
  \sum_{h=0}^{N-j}(-1)^h\binom{N-j}{h}\\
&=T^N.
\end{aligned}
\]
The identity $D_{n+1}(S,T)=S D_n(S,T)-D_n(S,T)T$
shows that $D_n(S,T)=0$ for all $n\geq p$. Hence
\[
T^N=\sum_{k=0}^{p-1}(-1)^k\binom Nk S^{N-k}D_k(S,T)
\qquad (N\geq p).
\]
Fix $x\in X$ and suppose that
$\alpha:=\limsup_n\|x\|_{T^n}^{1/n}<\infty$.
For any $c>\alpha$, choose $z_j\in X$ such that $T^jz_j=x$
and $\|z_j\|\leq c^j$ for all sufficiently large $j$.
For $n\geq1$, set
\[
y_n=\sum_{k=0}^{p-1}(-1)^k\binom{n+p-1}{k}
S^{p-1-k}D_k(S,T)z_{n+p-1}.
\]
Then $S^ny_n=x$. Using
$\binom{n+p-1}{k}\leq(n+p-1)^{p-1}$ for $0\leq k<p$, we obtain for all sufficiently large $n$
\[
\|y_n\|\leq c^{n+p-1}(n+p-1)^{p-1}
\sum_{k=0}^{p-1}\|S^{p-1-k}D_k(S,T)\|.
\]
Thus $\limsup_n\|x\|_{S^n}^{1/n}\leq c$, and letting
$c\downarrow\alpha$ gives the corresponding inequality.
The hypothesis $D_q(T,S)=0$ gives the reverse inequality
in the same way; it also covers the case where one of the
limits is infinite. Finally, suppose that $S=T+N$, where $N$ is nilpotent and
$N\in\operatorname{comm}_w(T)$. The commutator $ST-TS=NT-TN$ is annihilated on both sides by $T$ and $N$, hence also by $S$.
Thus the binomial formula of \cite[Corollary 3.6(i)]{GCP}, applied to $S,-T$ and to $T,-S$, gives
\[
D_n(S,T)=N^n,\qquad D_n(T,S)=(-N)^n\qquad(n\ge3).
\]
Both vanish for sufficiently large $n$.
\end{proof}

\begin{proposition}\label{onesidedpert}
Let $T,N\in\B(X)$ and put $S=T+N$.
\begin{enumerate}
\item If $N$ is nilpotent and $TNT=NT^2$, then $\mathcal I_r(T)\subseteq\mathcal I_r(S)$ for every $r\ge0$.
\item If $N$ is quasinilpotent and $TNT=NT^2$, then $\mathcal K_r(T)\subseteq\mathcal K_r(S)$ for every $r\ge0$.

\item If $N$ is quasinilpotent and $TNT=T^2N$, then
$\mathcal H_0(S)\subseteq\mathcal H_0(T)$.
\end{enumerate}
\end{proposition}
\begin{proof}
(i) Let $N^p=0$, $p\geq1$, and set $\ell=2(p-1)$.
Since $NT$ commutes with $T$, \cite[Lemma 3.1(II)(ii)]{GCP}
gives $(NT)^j=N^jT^j$ for $j\geq1$. Hence
\[
TN^jT^j=T(NT)^j=(NT)^jT=N^jT^{j+1}.
\]
Since $\ell-j\geq j$ for $0\leq j<p$, we obtain
\[
SN^jT^{\ell-j}=N^jT^{\ell-j+1}+N^{j+1}T^{\ell-j}
\quad(0\leq j<p).
\]
Define $Q_0=T^\ell$ and, for $n\geq1$, put
\[
Q_n=\sum_{j=0}^{p-1}(-1)^j
\binom{n+j-1}{j}N^jT^{\ell-j}.
\]
Pascal's identity gives $SQ_{n+1}=Q_nT$ for $n\geq0$.
Hence $S^nQ_n=T^{n+\ell}$. For $n\geq1$,
$\binom{n+j-1}{j}\leq(n+p-1)^{p-1}$ for $0\leq j<p$, so
\[
\|Q_n\|\leq(n+p-1)^{p-1}
\sum_{j=0}^{p-1}\|N^jT^{\ell-j}\|.
\]

Let $x\in\mathcal I_r(T)$ and choose $y_n$ such that
$T^ny_n=x$ with the required exponential growth bound.
Then $z_n=Q_ny_{n+\ell}$ satisfies $S^nz_n=x$ and
\[
\limsup_{n\to\infty}\|z_n\|^{1/n}
\leq\limsup_{n\to\infty}\|y_n\|^{1/n}.
\]
Thus $x\in\mathcal I_r(S)$, also for $r=0$.

(ii) Put $V=\mathcal R^\infty(T)$. Since $TNT=NT^2$, $T$ and $N$ commute on $\mathcal R(T)$, and $NT^{k+1}=T^kNT$ for $k\ge1$. Thus $NV\subseteq V$ and $TN^jv=N^jTv$ for $v\in V$.
Let $(y_n)$ be a backward chain for $T$ with $\|y_n\|\le M\theta^n$, where $\theta>0$ and $\theta<r^{-1}$ if $r>0$. Each $y_n$ belongs to $V$. Choose $\varepsilon>0$ with $\varepsilon\theta<1$ and, if $r>0$, $\theta/(1-\varepsilon\theta)<r^{-1}$. Since $N$ is quasinilpotent, $\|N^j\|\le C_\varepsilon\varepsilon^j$. Put $z_0=y_0$ and
\[
z_n=\sum_{j\ge0}(-1)^j\binom{n+j-1}{j}N^jy_{n+j}\qquad(n\ge1).
\]
The series converges in $X$ and the binomial series gives
$\|z_n\|\le MC_\varepsilon[\theta/(1-\varepsilon\theta)]^n$.
Applying $S=T+N$ termwise, the terms of positive order cancel for $Sz_1$, so $Sz_1=y_0$. For $n\ge1$, the coefficient of $(-1)^j N^jy_{n+j}$ in $Sz_{n+1}$ is $1$ when $j=0$ and
$\binom{n+j}{j}-\binom{n+j-1}{j-1}=\binom{n+j-1}{j}$ when $j\ge1$. Hence $Sz_{n+1}=z_n$ and $y_0\in\mathcal K_r(S)$.

(iii) Put $Q=TN$. The identity $TNT=T^2N$ implies
$QT=TQ$ and $Q^n=T^nN^n$ by \cite[Lemma 3.1(I)(ii)]{GCP},
so $Q$ is quasinilpotent.
For $A\in\B(X)$, write
$r(A,x)=\limsup_{n\to\infty}\|A^nx\|^{1/n}$.
Since $Q$ commutes with $T^2$ and is quasinilpotent,
\cite[Lemma 1(6), (7)]{Danes} gives
$r(T^2+Q,x)=r(T^2,x)=r(T,x)^2$.
As $TS=T^2+Q$ and $(TS)^n=T^nS^n$, we obtain
\[
r(T,x)^2=r(TS,x)
\leq\|T\|r(S,x).
\]
Therefore $r(S,x)=0$ implies $r(T,x)=0$.
\end{proof}

\smallskip
\noindent The following example shows that the inclusions in
Proposition \ref{onesidedpert} can be strict.
\begin{example}\label{onesidedstrict}
On $\C^3$, let $T=E_{13}$ and $N=E_{21}+E_{32}$, where
$E_{ij}$ are the matrix units. Then $T^2=N^3=TNT=0$, whereas
$(T+N)^3=I$. Thus, for $0\leq r<1$,
$\mathcal I_r(T)=\mathcal K_r(T)=\{0\}\subsetneq
\mathcal I_r(T+N)=\mathcal K_r(T+N)=\C^3$, and
$\mathcal H_0(T+N)=\{0\}\subsetneq\mathcal H_0(T)=\C^3$.
\end{example}

\section{Mixed algebraic decompositions}
\begin{theorem}\label{mixed}
Let $T\in\B(X)$, $a,b\geq0$, and let $W$ be either
$\mathcal K_b(T)$ or $\mathcal I_b(T)$. Put $L=\mathcal H_a(T)+W$. Then:
\begin{enumerate}
\item $L$ is either $X$ or meagre in $X$;
\item if $\codim L<\infty$, then $L$ is closed;
\item if $X=\mathcal H_a(T)\oplus W$ algebraically, both summands are closed and $W=\mathcal K_b(T)$.
\end{enumerate}
In (iii), $B=T_W$ is invertible. For $A=T_{\mathcal H_a(T)}$,
\begin{equation}\label{mixedbounds}
\begin{cases}
\rad(A)<a,&a>0,\\
\rad(A)=0,&a=0,
\end{cases}
\qquad
b>0\ \Longrightarrow\ \rad(B^{-1})<b^{-1}.
\end{equation}
\end{theorem}
\begin{proof}
Let $J=\mathbb N$ if $a>0$, and $J=\{1\}$ if $a=0$.
For $a>0$, enumerate $\mathbb Q\cap(0,a)$ as $(a_j)_{j\geq1}$,
put $Z_j=V_{a_j}$, and let $U_j:Z_j\to X$ be the inclusion.
For $a=0$, put $Z_1=\mathcal{H}_0(T)$ and let $U_1$ be the inclusion.
Next, enumerate $\mathbb Q\cap(0,b^{-1})$ as $(b_k)_{k\geq1}$
if $b>0$, and $\mathbb Q\cap(0,\infty)$ if $b=0$.
For $W=\mathcal K_b(T)$ set $(W_k,V_k)=(Y_{b_k},E_{b_k})$, and for $W=\mathcal I_b(T)$ set $(W_k,V_k)=(F_{b_k},J_{b_k})$. Propositions~\ref{models} and~\ref{indmodels} and \cite[Proposition 2.3]{AC} give continuous maps from Fr\'echet spaces with
\[
\mathcal H_a(T)=\bigcup_{j\in J}\mathcal R(U_j),\qquad W=\bigcup_{k\ge1}\mathcal R(V_k).
\]
For $j\in J$ and $k\geq1$, define
\[
\Phi_{j,k}:Z_j\times W_k\longrightarrow X,\qquad
\Phi_{j,k}(z,w)=U_jz+V_kw.
\]
Then $L=\bigcup_{j\in J,\ k\geq1}\mathcal{R}(\Phi_{j,k})$. If $L$ is nonmeagre, some $\mathcal{R}(\Phi_{j,k})$ is nonmeagre. By \cite[Theorem 2.11]{Rudin}, $\Phi_{j,k}$ is onto and open. Hence $L=X$.

Suppose $\codim L<\infty$. Choose a finite-dimensional algebraic complement $F$ of $L$ in $X$. For $j\in J$ and $k\geq1$, define
\[
\begin{split}
\Psi_{j,k}:Z_j\times W_k\times F&\longrightarrow X,\\
(z,w,f)&\longmapsto U_jz+V_kw+f.
\end{split}
\]
The product is a Fr\'echet space, and $\Psi_{j,k}$ is continuous. Since $X=\bigcup_{j\in J,\ k\geq1}\mathcal{R}(\Psi_{j,k})$, some $\Psi_{j,k}$ is an open surjection. Let $Q:X\to F$ be the algebraic projection onto $F$ along $L$. Then
\[
Q\Psi_{j,k}(z,w,f)=f
\qquad(z\in Z_j,\ w\in W_k,\ f\in F).
\]
Since $\Psi_{j,k}$ is open and onto, the continuity of $Q\Psi_{j,k}$ implies that $Q$ is continuous. Hence $L=\mathcal{N}(Q)$ is closed.

Suppose $X=\mathcal H_a(T)\oplus W$. Let $P$ be the algebraic projection onto $W$ along $\mathcal H_a(T)$. For an open surjection $\Phi_{j,k}$ as above, $P\Phi_{j,k}(z,w)=V_kw$. Thus $P$ is continuous and both summands are closed.

The two summands are invariant, so $PT=TP$. Since
$\mathcal N(T)\subseteq\mathcal H_a(T)$, $B=T_W$ is injective.
For $x\in W$, choose $x_n\in X$ with $T^nx_n=x$ and the
growth prescribed by $\mathcal I_b(T)$, since $W\subseteq\mathcal I_b(T)$. Then
\[
B^nPx_n=x\qquad(n\geq0).
\]
Taking $n=1$ shows that $B$ is surjective. Hence
\[
B^{-n}x=Px_n,\qquad
\|B^{-n}x\|\leq\|P\|\,\|x_n\|.
\]
Thus $x\in\mathcal K_b(T)$, so $W=\mathcal K_b(T)$.
On the inner summand the forward growth of each vector is less than $a$ when $a>0$, and zero when $a=0$. On $W$, the unique internal chain is $(B^{-n}x)$ and its growth is less than $b^{-1}$ when $b>0$. For a bounded operator $D$ on a nonzero Banach space $Y$, Dane\v{s}'s lemma gives $\rad(D)=\max_{y\in Y}\limsup_n\|D^ny\|^{1/n}$ \cite[Lemma 2]{Danes}. Apply this formula to $A$ and, for $b>0$, to $B^{-1}$ to obtain \eqref{mixedbounds}.
\end{proof}

\begin{corollary}\label{indcriterion}
Fix $a,b\ge0$ and let $W_b$ denote either $\mathcal K_b$ or $\mathcal I_b$. Then $X=\mathcal H_a(T)\oplus W_b(T)$ algebraically if and only if there is a reducing decomposition $T=A\oplus B$ on $X=M\oplus N$ such that
\[
\mathcal H_a(B)=W_b(A)=\{0\},\qquad B\text{ is invertible},
\]
$\rad(A)<a$ for $a>0$, $\rad(A)=0$ for $a=0$, and $\rad(B^{-1})<b^{-1}$ for $b>0$.
In this case $M=\mathcal H_a(T)$ and $N=W_b(T)=\mathcal K_b(T)$.
In particular, $X=\mathcal H_a(T)\oplus\mathcal I_b(T)$ holds exactly when
$X=\mathcal H_a(T)\oplus\mathcal K_b(T)$ and $\mathcal H_a(T)\cap\mathcal I_b(T)=\{0\}$.
\end{corollary}
\begin{proof}
For the forward implication use Theorem~\ref{mixed} and take the canonical summands. The two vanishing conditions follow by directness. Conversely, the radius bound on $A$ gives $\mathcal H_a(A)=M$, and $(B^{-n}x)$ gives $W_b(B)=N$. The bounded coordinate projections commute with $T$, so both radial spaces split over $M\oplus N$. The last assertion follows from Theorem~\ref{mixed} and $\mathcal K_b(T)\subseteq\mathcal I_b(T)$.
\end{proof}

\begin{corollary}\label{canonical}
Let $0<r\leq s$. The following assertions are equivalent:
\begin{enumerate}
\item $X=\mathcal{H}_r(T)\oplus \mathcal{K}_s(T)$ algebraically;
\item $\sigma(T)\cap\{\lambda:r\leq|\lambda|\leq s\}=\varnothing$;
\item $X=\mathcal{H}_r(T)\oplus \mathcal{K}_s(T)$ topologically.
\end{enumerate}
Whenever these conditions hold, for every $t\in[r,s]$,
\begin{equation}\label{radiuscollapse}
\mathcal{H}_t(T)=\mathcal{H}_r(T),\qquad \mathcal{K}_t(T)=\mathcal{K}_s(T),\qquad
X=\mathcal{H}_t(T)\oplus \mathcal{K}_t(T).
\end{equation}
For $s=r$, these conditions are also equivalent to
\[
\mathcal{H}_r(T)+\mathcal{K}_r(T)\text{ is nonmeagre},\qquad
\dim(\mathcal{H}_r(T)\cap \mathcal{K}_r(T))<\infty.
\]
\end{corollary}
\begin{proof}
If (i) holds, Theorem \ref{mixed} gives
\[
\rad(T_{\mathcal{H}_r(T)})<r,\qquad
\rad((T_{\mathcal{K}_s(T)})^{-1})<s^{-1}.
\]
Hence (ii) and (iii) follow. Conversely, assume (ii). Since the annulus
$\{\lambda:r\leq|\lambda|\leq s\}$ is disjoint from $\sigma(T)$,
the inner Riesz projection is independent of $t\in[r,s]$.
On the inner block $A$, $\rad(A)<t$ and $\mathcal K_t(A)\subseteq\mathcal I_t(A)=\{0\}$ by Remark~\ref{indradremark}. On the outer block $B$, $\rad(B^{-1})<t^{-1}$ gives $\mathcal K_t(B)$ equal to the whole block, while $\|x\|\le\|B^{-n}\|\,\|B^nx\|$ gives $\mathcal H_t(B)=\{0\}$. Thus the range and kernel of the projection are $\mathcal H_t(T)$ and $\mathcal K_t(T)$, respectively.
This gives \eqref{radiuscollapse}, and (i) follows.
Finally, for $s=r$, apply Theorem \ref{mixed}(i) and Proposition
\ref{dichotomy}.
\end{proof}

\begin{corollary}\label{annuli}
Suppose $0<s<r$ and $X=\mathcal{H}_t(T)\oplus \mathcal{K}_t(T)$ algebraically for $t=s,r$. For $t=s,r$, denote by $Q_t$ the inner Riesz projection of $T$ at radius $t$ (the circles lie in $\rho(T)$ by Corollary~\ref{canonical}). Then
\begin{align}\label{annularsum}
\mathcal{K}_s(T)&=\mathcal{K}_r(T)\oplus F_{s,r},&
\mathcal{H}_r(T)&=\mathcal{H}_s(T)\oplus F_{s,r},\\
F_{s,r}&=\mathcal{K}_s(T)\cap \mathcal{H}_r(T)=\mathcal{R}(Q_r-Q_s).\notag
\end{align}
All these subspaces are closed. The restriction $T_{F_{s,r}}$ is invertible and, for $F_{s,r}\ne\{0\}$,
\[
\sigma(T_{F_{s,r}})=\sigma(T)\cap\{\lambda:s<|\lambda|<r\}.
\]
The quotient operator on $\mathcal{K}_s(T)/\mathcal{K}_r(T)$ is similar to $T_{F_{s,r}}$.
\end{corollary}
\begin{proof}
By Corollary~\ref{canonical}, $\mathcal R(Q_t)=\mathcal H_t(T)$ and $\mathcal N(Q_t)=\mathcal K_t(T)$ for $t=s,r$.

Since $Q_sQ_r=Q_rQ_s=Q_s$, the operator $Q_r-Q_s$ is a
bounded projection with range $F_{s,r}$, and
\eqref{annularsum} follows. The quotient identification is
$f\mapsto f+\mathcal K_r(T)$ for $f\in F_{s,r}$.
\end{proof}
\begin{corollary}\label{H0Ks}
Let $s>0$. The following assertions are equivalent:
\begin{enumerate}
\item $X=\mathcal{H}_0(T)\oplus \mathcal{K}_s(T)$ algebraically;
\item $\sigma(T)\cap\{\lambda:0<|\lambda|\leq s\}=\varnothing$;
\item $T$ is generalized Drazin invertible and $\mathcal{K}_s(T)=\mathcal{K}(T)$.
\end{enumerate}
These conditions imply that the sum is topological and that, for every $0<t\leq s$,
\[
\mathcal{H}_t(T)=\mathcal{H}_0(T),\qquad \mathcal{K}_t(T)=\mathcal{K}(T).
\]
\end{corollary}
\begin{proof}
Under (i), Theorem \ref{mixed} yields $T=A\oplus B$, with $\rad(A)=0$, $B$ invertible and $\rad(B^{-1})<s^{-1}$, giving (ii). Under (ii), the spectral decomposition at $0$ gives $\mathcal{H}_0(T)$ and $\mathcal{K}(T)$ as its inner and outer spaces, respectively \cite[Theorem 7.1]{Koliha} and \cite[Theorem 3.1]{AC}; the circular projections are constant for $0<t\leq s$. Thus (iii) follows. The implication (iii)$\Rightarrow$(i) follows from \cite[Theorem 3.1]{AC}.
\end{proof}

\smallskip
\noindent Unlike the case $0\leq r<s$, the decomposition
$X=\mathcal H_r(T)\oplus\mathcal K_s(T)$ for $0\leq s<r$
cannot be characterized by the spectrum alone, as the following
example shows.

\begin{example}\label{shifts}
Let $0\leq s<c<r$ and let $S$ be the unilateral right shift on $\ell^2(\mathbb N)$. The operators $cS$ and $cS^*$ have the same spectrum, the closed disc of radius $c$, and
\[
\mathcal{H}_r(cS)=\mathcal{H}_r(cS^*)=\ell^2(\mathbb N).
\]
However, $\mathcal{K}_s(cS)=\{0\}$, since $\mathcal{K}_s(cS)\subset\bigcap_{n\geq1}\mathcal{R}(S^n)=\{0\}$, whereas $\mathcal{K}_s(cS^*)=\ell^2(\mathbb N)$: for each $x$, the chain $x_n=c^{-n}S^nx$ has growth $c^{-1}<s^{-1}$ if $s>0$, and is exponentially bounded if $s=0$. Hence $X=\mathcal{H}_r(cS)\oplus \mathcal{K}_s(cS)$, but the analogous identity fails for $cS^*$. 
\end{example}

\section{Decompositions with independent preimages}\label{sec:indecomp}

\begin{proposition}\label{indordered}
Let $0\leq a\leq b$. Then
\[
X=\mathcal H_a(T)\oplus\mathcal I_b(T)
\quad\text{algebraically}
\iff
X=\mathcal H_a(T)\oplus\mathcal K_b(T)
\quad\text{algebraically}.
\]
Whenever these conditions hold, the sum is topological and
$\mathcal I_b(T)=\mathcal K_b(T)$. In particular:
\begin{enumerate}
\item for $0<r\leq s$, the decomposition
$X=\mathcal H_r(T)\oplus\mathcal I_s(T)$ holds if and only if
$\sigma(T)\cap\{\lambda:r\leq|\lambda|\leq s\}=\varnothing$;
\item for $s>0$, the decomposition
$X=\mathcal H_0(T)\oplus\mathcal I_s(T)$ holds if and only if
$\sigma(T)\cap\{\lambda:0<|\lambda|\leq s\}=\varnothing$;
\item $X=\mathcal H_0(T)\oplus\mathcal I(T)$ holds if and only if
$T$ is generalized Drazin invertible.
\end{enumerate}
In each statement, an algebraic direct sum may equivalently be
required to be topological.
\end{proposition}
\begin{proof}
The forward implication follows from Theorem~\ref{mixed}. For the converse, let $X=\mathcal H_a(T)\oplus\mathcal K_b(T)$ and write $A=T_{\mathcal H_a(T)}$. Theorem~\ref{mixed} gives $\rad(A)\le b$. Hence $\mathcal I_b(A)=\{0\}$ by Remark~\ref{indradremark}, and Corollary~\ref{indcriterion} applies. The three spectral criteria follow from Corollaries~\ref{canonical} and~\ref{H0Ks} and \cite[Theorem 3.1]{AC}.
\end{proof}

\smallskip
\noindent For $b<a$, a decomposition with $\mathcal K_b(T)$ need not remain direct after replacing $\mathcal K_b(T)$ by $\mathcal I_b(T)$.
Indeed, for the operator in Example \ref{indcounterexample}, choose
$a>\max\{1,\|T\|\}$ and $b=1$. Then
$\mathcal H_a(T)=X$ and $\mathcal K_b(T)=\{0\}$, whereas
$\mathcal I_b(T)=\C e_0\ne\{0\}$.

\begin{theorem}\label{closedindcores}
Let $X$ be a complex Banach space, $T\in\B(X)$, and
$0\leq a\leq b$. If $\mathcal I_b(T)$ and
$\mathcal H_a(T)+\mathcal I_b(T)$ are closed, then
\[
\mathcal I_b(T)=\mathcal K_b(T).
\]
\end{theorem}
\begin{proof}
Put $E=\mathcal I_b(T)$ and $L=\mathcal H_a(T)+E$.
We may assume that $E\ne\{0\}$. Since $E$ is closed, Proposition \ref{indmodels} expresses $E$
as a countable union of ranges of the maps $J_\eta:F_\eta\to E$.
The nonmeagre-range
open mapping theorem therefore gives $C\geq1$ and $\eta>0$,
with $\eta<b^{-1}$ if $b>0$, such that for every $x\in E$
and $n\geq0$ there is $y\in X$ satisfying
\begin{equation}\label{controlledindpreimages}
T^ny=x,\qquad \|y\|\leq C\eta^n\|x\|.
\end{equation}
Suppose $a>0$. Since $L$ is closed, Proposition \ref{models} gives
\[
L=\bigcup_{\substack{t\in\mathbb Q\\0<t<a}}(V_t+E).
\]
The same theorem, applied to the maps
$V_t\times E\to L$, $(h,i)\mapsto h+i$, gives $0<t<a$
and $D>0$ such that every $y\in L$ can be written
$y=h+i$, with $h\in V_t$, $i\in E$, and
\begin{equation}\label{decompcontrol}
\|i\|+\sup_{j\geq0}t^{-j}\|T^jh\|\leq D\|y\|.
\end{equation}
When $a=0$, the open mapping theorem applied to the surjection
$\mathcal H_0(T)\times E\to L$ yields the same estimate for every
$t>0$.

Let $S$ be induced by $T$ on $X/E$. We first prove
$\mathcal K_b(S)=\{0\}$. Take a backward chain $(z_n)$ for $S$
with $\|z_n\|\leq M\theta^n$, where $0<\theta<b^{-1}$ if
$b>0$, and $\theta$ is finite and exceeds $\limsup_n\|z_n\|^{1/n}$ if $b=0$. Choose representatives
$x_n\in X$ with $\|x_n\|\leq2M\theta^n$. Then
\[
e_j=Tx_{j+1}-x_j\in E,\qquad
\|e_j\|\leq M_1\theta^j,\qquad
M_1=2M(1+\|T\|\theta).
\]
By \eqref{controlledindpreimages}, for $0\leq j<n$ choose
$w_{n,j}\in X$ such that
\[
T^{n-j}w_{n,j}=e_j,\qquad
\|w_{n,j}\|\leq C\eta^{n-j}\|e_j\|.
\]
Since
\[
T^nx_n-x_0=\sum_{j=0}^{n-1}T^je_j,
\]
the vector $v_n=x_n-\sum_{j=0}^{n-1}w_{n,j}$ satisfies
\[
T^nv_n=x_0,\qquad
\|v_n\|\leq2M\theta^n+
CM_1\sum_{j=0}^{n-1}\eta^{n-j}\theta^j
\leq(2M+CM_1n)q^n,
\]
where $q=\max\{\eta,\theta\}$. Thus
\[
\limsup_{n\to\infty}\|v_n\|^{1/n}
\leq\max\{\eta,\theta\},
\]
which is less than $b^{-1}$ if $b>0$, and finite if $b=0$.
Thus $x_0\in\mathcal I_b(T)=E$, so $z_0=0$.

Let $\pi:X\to X/E$ be the quotient map and put $Z=L/E$. This is a closed $S$-invariant subspace. If $a>0$, take $t$ from \eqref{decompcontrol}; if $a=0$, choose $0<t<\eta^{-1}$. In either case $t\eta<1$. For each $z\in Z$, take a representative $y\in L$ with $\|y\|\le2\|z\|$. Writing $y=h+i$ as in \eqref{decompcontrol} gives
\[
\|S^jz\|=\|\pi T^jh\|\le 2Dt^j\|z\|\qquad(j\ge0).
\]
For $x\in E$ choose $y_n$ as in \eqref{controlledindpreimages}, with $y_0=x$. Set $p_n=\pi y_n$ and $f_n=Sp_{n+1}-p_n$. Since $Ty_{n+1}-y_n\in\mathcal N(T^n)\subseteq L$, we have $f_n\in Z$ and $\|f_n\|\le C(1+\|T\|\eta)\eta^n\|x\|$. The series
\[
u_n=\sum_{j\ge0} S^jf_{n+j}
\]
converges in $Z$, with $\|u_n\|\le 2DC(1+\|T\|\eta)(1-t\eta)^{-1}\eta^n\|x\|$. As $u_n=f_n+Su_{n+1}$, the sequence $(p_n+u_n)$ is a backward chain for $S$. Every tail has growth at most $\eta$. Thus $p_n+u_n\in\mathcal K_b(S)=\{0\}$ for every $n$, so $\pi y_n=-u_n\in Z$ and $y_n\in L$.

Write $y_n=h_n+i_n$ as in \eqref{decompcontrol}. We obtain

\[
\|i_n\|\leq CD\eta^n\|x\|,\qquad
\|x-T^ni_n\|=\|T^nh_n\|
\leq CD(t\eta)^n\|x\|.
\]
Choose an integer $m\geq1$ such that
\[
\delta=CD(t\eta)^m<1,\qquad
d=\frac{CD\eta^m}{1-\delta},
\]
and, when $b>0$, $d^{1/m}<b^{-1}$. Such a choice is possible
since $t\eta<1$ and $d^{1/m}\longrightarrow\eta$.
Therefore, for every $x\in E$, there exists $i\in E$ such that
\[
\|x-T^mi\|\leq\delta\|x\|,\qquad
\|i\|\leq CD\eta^m\|x\|.
\]
The successive-approximation argument in the proof of
Theorem \ref{indsurjspectrum}, applied to $T_E^m$, gives, for every
$x\in E$, an element $u\in E$ such that
\[
T^mu=x,\qquad\|u\|\leq d\|x\|.
\]
Iterate this estimate under $T_E^m$ and fill the $m-1$ intermediate indices with forward powers of $T$, as in Proposition~\ref{indmodels}. The resulting backward chain has growth at most $d^{1/m}$.
Hence $x\in\mathcal K_b(T)$, and
$\mathcal I_b(T)=\mathcal K_b(T)$.
\end{proof}

\smallskip
\noindent Recall that an operator $S\in\B(X)$ is semi-regular if
$\mathcal R(S)$ is closed and
$\mathcal N(S)\subseteq\mathcal R^\infty(S)$.
An operator $T\in\B(X)$ is pseudo-Fredholm if it admits a
reducing decomposition $T=S\oplus Q$, where $S$ is semi-regular
and $Q$ is quasinilpotent. If $S$ can be chosen bounded below (respectively, surjective),
then $T$ is called left (respectively, right) generalized Drazin invertible.

\begin{corollary}\label{indpseudofredholm}
Let $T\in\B(X)$. The following assertions hold.
\begin{enumerate}
\item If $\mathcal I(T)$ and
$\mathcal I(T)+\mathcal H_0(T)$ are closed, and both
$\mathcal N(T)\cap\mathcal I(T)$ and
$\mathcal R(T)+\mathcal H_0(T)$ are complemented in $X$,
then $T$ is pseudo-Fredholm.
\item If $X=H$ is a complex Hilbert space, then $T$ is
pseudo-Fredholm if and only if $\mathcal I(T)$,
$\mathcal I(T)+\mathcal H_0(T)$ and
$\mathcal R(T)+\mathcal H_0(T)$ are closed.
\end{enumerate}
Every pseudo-Fredholm operator on a complex Banach space satisfies $\mathcal I(T)=\mathcal K(T)$.
\end{corollary}
\begin{proof}
Under the first two closedness assumptions, Theorem \ref{closedindcores}
with $a=b=0$ gives $\mathcal I(T)=\mathcal K(T)$. The Banach-space
assertion follows from \cite[Theorem 3.1]{ACPF}.

On a Hilbert space, the closedness conditions in (ii) imply
pseudo-Fredholmness, since closed subspaces are complemented.
For any pseudo-Fredholm operator on a Banach space, write $T=S\oplus Q$ with $S$ semi-regular and $Q$ quasinilpotent.
Since $\mathcal K(S)=\mathcal R^\infty(S)$ and
$\mathcal I(S)\subseteq\mathcal R^\infty(S)$, we have
$\mathcal I(S)=\mathcal K(S)$. Also,
$\mathcal I(Q)=\mathcal K(Q)=\{0\}$ by Remark \ref{indradremark}.
Hence $\mathcal I(T)=\mathcal K(T)$, and (ii) follows from
\cite[Corollary 3.2]{ACPF}.
\end{proof}

\begin{remark}
The three closedness conditions in (ii) are independent. Indeed,
the operators in \cite[Proposition 6.3]{ACPF} are injective, so that
$\mathcal I(T)=\mathcal K(T)$ for each of them.
\end{remark}

\begin{corollary}\label{indleftgd}
Let $T\in\B(X)$. The following assertions hold.
\begin{enumerate}
\item If $\mathcal I(T)$ and
$\mathcal I(T)+\mathcal H_0(T)$ are closed and
$\mathcal N(T)\cap\mathcal I(T)=\{0\}$,
then $\mathcal H_0(T)$ is closed. If, in addition,
$\mathcal R(T)+\mathcal H_0(T)$ is complemented in $X$,
then $T$ is left generalized Drazin invertible.
\item If $\mathcal I(T)$ is closed,
$X=\mathcal I(T)+\mathcal H_0(T)$ and
$\mathcal N(T)\cap\mathcal I(T)$ is complemented in $X$,
then $T$ is right generalized Drazin invertible.
\end{enumerate}
If $X=H$ is a complex Hilbert space, then:
\begin{enumerate}[start=3]
\item $T$ is left generalized Drazin invertible if and only if
$\mathcal H_0(T)$ and $\mathcal R(T)+\mathcal H_0(T)$ are closed.
Equivalently, $\mathcal I(T)$, $\mathcal I(T)+\mathcal H_0(T)$ and
$\mathcal R(T)+\mathcal H_0(T)$ are closed, and
$\mathcal N(T)\cap\mathcal I(T)=\{0\}$.
\item $T$ is right generalized Drazin invertible if and only if
$\mathcal I(T)$ is closed and
$H=\mathcal I(T)+\mathcal H_0(T)$.
\end{enumerate}
In each case, $\mathcal I(T)=\mathcal K(T)$.
\end{corollary}
\begin{proof}
Under the first assumptions in (i), put
$E=\mathcal I(T)$ and $L=E+\mathcal H_0(T)$.
By Theorem \ref{closedindcores}, $E=\mathcal K(T)$.
Since $T(E)=E$ and $\mathcal N(T)\cap E=\{0\}$,
$T_E$ is invertible. Hence
$E\cap\mathcal H_0(T)=\mathcal H_0(T_E)=\{0\}$.
The addition map
\[
E\times\mathcal H_0(T)\longrightarrow L,\qquad
(e,h)\longmapsto e+h
\]
is therefore a continuous linear bijection between Fr\'echet
spaces. By the open mapping theorem it is a topological
isomorphism. Thus $\mathcal H_0(T)$ is closed in $L$, and hence
in $X$.

Under the Banach-space assumptions, Corollary \ref{indpseudofredholm}
gives a generalized Kato decomposition $T=S\oplus Q$ with
$\mathcal I(T)=\mathcal K(T)$. Since
$\mathcal N(S)\subseteq\mathcal K(S)$, the kernel condition implies
$\mathcal N(S)=\{0\}$. Thus $S$ is bounded below.

Under the assumptions in (ii),
Theorem \ref{closedindcores} gives
$\mathcal I(T)=\mathcal K(T)$, while
$X=\mathcal I(T)+\mathcal H_0(T)$ implies
$X=\mathcal R(T)+\mathcal H_0(T)$.
By Corollary \ref{indpseudofredholm}, $T$ admits a generalized
Kato decomposition $T=S\oplus Q$ on $X=M\oplus N$.
Since $\mathcal H_0(S)\subseteq\mathcal K(S)$,
we have $M=\mathcal K(S)+\mathcal H_0(S)=\mathcal K(S)$.
Thus $S$ is surjective.

On a Hilbert space, the first equivalence in (iii) is
\cite[Lemma 2.3]{ACZReg}. The conditions involving $\mathcal I(T)$ imply left generalized Drazin invertibility by (i). Conversely, in a decomposition $T=D\oplus Q$ with $D$ bounded below and $Q$ quasinilpotent, Corollary~\ref{indpseudofredholm} gives the three closedness conditions. Since $\mathcal N(D)=\{0\}$ and $\mathcal I(Q)=\{0\}$, we also have $\mathcal N(T)\cap\mathcal I(T)=\{0\}$.

For (iv), a reducing decomposition into a surjective operator
and a quasinilpotent operator gives
$\mathcal I(T)=\mathcal K(T)$ closed and
$H=\mathcal I(T)+\mathcal H_0(T)$.
The converse follows from the right generalized Drazin
assertion in (ii), since closed subspaces of $H$ are complemented.
\end{proof}

\section{Finite defects and annular decompositions}

\begin{lemma}\label{saturation}
Let $T\in\B(X)$ and $r,s\geq0$. For
$W_s(T)\in\{\mathcal K_s(T),\mathcal I_s(T)\}$, put
$L=\mathcal H_r(T)+W_s(T)$. Then
\[
(T-\lambda I)^{-1}(L)=L
\qquad\text{if }|\lambda|\geq r\text{ and }|\lambda|>s.
\]
For $W_s(T)=\mathcal K_s(T)$, this equality also holds when
$\lambda=0$ or $0<|\lambda|<\min\{r,s\}$.
If $\codim L<\infty$, then $L$ is closed, and the operator
$\widetilde T$ induced by $T$ on $X/L$ satisfies
\[
\sigma(\widetilde T)\subseteq
\begin{cases}
\{\lambda:\min\{r,s\}\leq|\lambda|\leq\max\{r,s\}\}
   \setminus\{0\},&W_s(T)=\mathcal K_s(T),\\
\{\lambda:|\lambda|<r\text{ or }|\lambda|\leq s\},
   &W_s(T)=\mathcal I_s(T).
\end{cases}
\]
\end{lemma}
\begin{proof}
Write $H=\mathcal H_r(T)$ and $W=W_s(T)$.
Both spaces are $T$-invariant, so $L\subseteq(T-\lambda I)^{-1}(L)$.

Suppose $|\lambda|\geq r$ and $|\lambda|>s$, and let
$(T-\lambda I)x=h+k$ with $h\in H$ and $k\in W$.
For $r>0$, choose
$\limsup_n\|T^nh\|^{1/n}<q<r\leq|\lambda|$.
When $r=0$, take any $0<q<|\lambda|$.
Since $\|T^nh\|\leq Cq^n$ for some $C>0$, the series
\[
u=-\sum_{j\geq0}\lambda^{-j-1}T^jh
\]
converges and satisfies $(T-\lambda I)u=h$.
Moreover,
$\|T^nu\|\leq Cq^n/(|\lambda|-q)$.
Hence $u\in H$; for $r=0$ this follows by taking $q$ arbitrarily small.
Put $z=x-u$, so that $(T-\lambda I)z=k$.

Choose $k_0=k$ and $k_n\in X$ such that $T^nk_n=k$ and
\[
\alpha:=\limsup_{n\to\infty}\|k_n\|^{1/n}<s^{-1}
\]
if $s>0$, or $\alpha<\infty$ if $s=0$.
If $W=\mathcal K_s(T)$, take a backward chain
$Tk_{n+1}=k_n$. For $n\geq1$, put
\[
z_n=\lambda^{-n}z-
\sum_{j=1}^n\lambda^{-(n-j+1)}k_j,
\qquad z_0=z.
\]
From $Tz=\lambda z+k$ and $T^jk_j=k$, we obtain
$T^nz_n=z$. If $(k_n)$ is a backward chain, then
$Tz_{n+1}=z_n$ as well.

Since $|\lambda|>s$, we may choose
$q>\max\{|\lambda|^{-1},\alpha\}$ with $q<s^{-1}$
when $s>0$, and a finite such $q$ when $s=0$.
Then $\|k_n\|\leq Cq^n$, and
\[
\|z_n\|\leq|\lambda|^{-n}\|z\|+
\frac{Cq^n}{|\lambda|}\sum_{j=0}^{n-1}(|\lambda|q)^{-j}
\leq C'q^n,
\]
since $|\lambda|q>1$. Thus $z\in W$, and $x\in L$.

Now let $W=\mathcal K_s(T)$. If $\lambda=0$ and
$Tx=h+k\in L$, choose $v\in W$ with $Tv=k$, since
$T\mathcal K_s(T)=\mathcal K_s(T)$. Then
$T(x-v)=h\in H$, and $T^{-1}(H)=H$ gives $x-v\in H$.
Hence $x\in L$.

Suppose $0<|\lambda|<\min\{r,s\}$, and write
$(T-\lambda I)x=h+k$ with $h\in H$ and $k\in W$.
Choose a chain $(k_n)$ for $k$ with $\|k_n\|\le C\theta^n$ and $\theta<s^{-1}$. Since $|\lambda|\theta<1$, the vectors $v_n=\sum_{j\ge0}\lambda^j k_{n+j+1}$ form a backward chain with $\|v_n\|\le C\theta^{n+1}/(1-|\lambda|\theta)$. Thus $v=v_0\in W$ and $(T-\lambda I)v=k$. Set
$w=x-v$. Then $(T-\lambda I)w=h$, and
\[
T^nw=\lambda^nw+\sum_{j=0}^{n-1}
\lambda^{n-1-j}T^jh.
\]
Choose $q$ with
$\max\{|\lambda|,\limsup_n\|T^nh\|^{1/n}\}<q<r$.
The estimate $\|T^jh\|\leq Cq^j$ now gives
$\|T^nw\|\leq C'q^n$, so $w\in H$ and $x\in L$.

Finally, if $\codim L<\infty$, then $L$ is closed by
Theorem \ref{mixed}(ii). For each value of $\lambda$ considered above,
\[
\mathcal N(\widetilde T-\lambda I)
=\{x+L:(T-\lambda I)x\in L\}=\{0\}.
\]
Thus $\widetilde T-\lambda I$ is invertible on the
finite-dimensional space $X/L$. The spectral inclusions follow.
\end{proof}

\begin{lemma}\label{indfiniterigidity}
Let $T\in\B(X)$ and $a,b\geq0$. Let
$W=\mathcal K_b(T)$, or $W=\mathcal I_b(T)$ when $a\leq b$.
Suppose that
\[
\dim(\mathcal H_a(T)\cap W)<\infty,\qquad
\codim(\mathcal H_a(T)+W)<\infty.
\]
Then $\mathcal H_a(T)$ and $W$ are closed, $T_W$ is invertible, and
there exists $\varepsilon>0$ such that
\begin{equation}\label{localmixed}
\sigma_*(T)\cap\{0<|\lambda|<\varepsilon\}
=\sigma_*(T_{\mathcal H_a(T)})\cap\{0<|\lambda|<\varepsilon\}
\end{equation}
for $\sigma_*\in\{\sigma,\sigma_b,\sigma_d\}$.
If $W=\mathcal I_b(T)$, then
\[
\mathcal H_a(T)\cap W=\{0\},\qquad W=\mathcal K_b(T).
\]
If $a=b=0$, then $\mathcal H_0(T)+W=X$.
\end{lemma}
\begin{proof}
Put $H=\mathcal H_a(T)$, $F=H\cap W$, and $L=H+W$.
By Theorem \ref{mixed}(ii), $L$ is closed. The range
representations in Propositions \ref{models} and \ref{indmodels},
with the Fr\'echet topology of $\mathcal H_0(T)$ when $a=0$,
give an open surjection
\[
\Phi:Z\times Y\to L,\qquad \Phi(z,y)=Uz+Vy,
\]
where $Z,Y$ are Fr\'echet spaces and $U,V$ are continuous,
with $\mathcal R(U)\subseteq H$ and $\mathcal R(V)\subseteq W$.
Since $F$ is finite-dimensional, the maps
\[
P_H(h+w)=h+F,\qquad P_W(h+w)=w+F
\]
from $L$ to $X/F$ are well defined. Their compositions with
$\Phi$ are continuous, and $\Phi$ is open and onto.
Thus $P_H$ and $P_W$ are continuous. As
$\mathcal N(P_H)=W$ and $\mathcal N(P_W)=H$, both $H$ and $W$
are closed.

Suppose that $W=\mathcal I_b(T)$, with $a\leq b$.
Since $W$ and $L$ are closed, Theorem \ref{closedindcores}
gives $W=\mathcal K_b(T)$. Proposition \ref{dichotomy} then
yields $F=\{0\}$. Thus $W=\mathcal K_b(T)$ in either case.

By Proposition \ref{dichotomy}, $T_F$ is invertible when
$F\ne\{0\}$, and trivially so when $F=\{0\}$. Since
$\mathcal N(T_W)\subseteq F$ and $TW=W$, the operator $T_W$
is invertible.

Put $K=W$ and $S=T_H\oplus T_K$ on $H\times K$. The
bounded surjection
\[
J:H\times K\longrightarrow L,\qquad J(h,k)=h+k
\]
satisfies $JS=T_LJ$ and
\[
G:=\mathcal N(J)=\{(f,-f):f\in F\}.
\]
Since $G$ is finite-dimensional and $S$-invariant, a
closed complement of $G$ in $H\times K$ gives
\[
S=\begin{pmatrix}S_G&B_1\\0&A\end{pmatrix},
\qquad S_G\simeq T_F,\quad A\simeq T_L.
\]
Similarly, $L$ is closed and of finite codimension. Hence
$X=L\oplus Z$ for some finite-dimensional $Z$, and
\[
T=\begin{pmatrix}T_L&B_2\\0&D\end{pmatrix},
\qquad D\simeq\widetilde T,
\]
where $\widetilde T$ is induced by $T$ on $X/L$.
Since $G$ and $Z$ are finite-dimensional, these triangular
representations give
\[
\begin{aligned}
\sigma(S)&=\sigma(S_G)\cup\sigma(T_L),&
\sigma(T)&=\sigma(T_L)\cup\sigma(D),\\
\sigma_b(S)&=\sigma_b(T_L)=\sigma_b(T),&
\sigma_d(S)&=\sigma_d(T_L)=\sigma_d(T).
\end{aligned}
\]
Since $S=T_H\oplus T_K$ on $H\times K$,
\[
\sigma_*(S)=\sigma_*(T_H)\cup\sigma_*(T_K),
\qquad \sigma_*\in\{\sigma,\sigma_b,\sigma_d\}.
\]
As $T_K$ is invertible and $\sigma(S_G)\cup\sigma(D)$
is finite, choose $\varepsilon>0$ such that
\[
\bigl(\sigma(S_G)\cup\sigma(D)\cup\sigma(T_K)\bigr)
\cap U=\varnothing,\qquad
U=\{\lambda:0<|\lambda|<\varepsilon\}.
\]
Consequently,
\[
\sigma_*(T)\cap U
=\sigma_*(T_L)\cap U
=\sigma_*(S)\cap U
=\sigma_*(T_H)\cap U,
\qquad \sigma_*\in\{\sigma,\sigma_b,\sigma_d\}.
\]
Finally, if $a=b=0$, we have $W=\mathcal K(T)$.
Lemma \ref{saturation} gives
$(T-\lambda I)^{-1}(L)=L$ for every $\lambda\in\C$.
Thus $\widetilde T-\lambda I$ is injective for every
$\lambda\in\C$. Since $X/L$ is finite-dimensional, it
follows that $L=X$.
\end{proof}

\begin{theorem}\label{finitecircle}
Let $T\in\B(X)$ and $r>0$. The following assertions are equivalent:
\begin{enumerate}
\item $\dim(\mathcal{H}_r(T)\cap \mathcal{K}_r(T))<\infty$ and
$\codim(\mathcal{H}_r(T)+\mathcal{K}_r(T))<\infty$;
\item $\mathbb T_r\cap\sigma_b(T)=\varnothing$, or equivalently, the spectral points on $\mathbb T_r$ form a finite set of isolated points with finite-rank Riesz projections;
\item there is a finite-dimensional $T$-invariant subspace $F_r$ such that
\[
X=\mathcal{H}_r(T)\oplus F_r\oplus \mathcal{K}_r(T)
\]
topologically and $\sigma(T_{F_r})\subset\mathbb T_r$.
\item $\dim(\mathcal H_r(T)\cap\mathcal I_r(T))<\infty$ and
$\codim(\mathcal H_r(T)+\mathcal I_r(T))<\infty$.
\end{enumerate}
The empty spectral set and $F_r=\{0\}$ are allowed. Under these conditions the intersection is zero, both radial spaces are closed, $\mathcal I_r(T)=\mathcal K_r(T)$, and
\begin{equation}\label{defectrank}
\codim(\mathcal{H}_r(T)+\mathcal{K}_r(T))=\dim F_r
=\sum_{\lambda\in\sigma(T)\cap\mathbb T_r}\operatorname{rank}P_\lambda.
\end{equation}
Here $F_r$ is the spectral subspace associated with $\sigma(T)\cap\mathbb T_r$, and $P_\lambda$ is the Riesz projection of $T$ at the isolated spectral point $\lambda$.
\end{theorem}
\begin{proof}
The two descriptions in (ii) agree by compactness of $\sigma(T)\cap\mathbb T_r$ and the definition of $\sigma_b(T)$.
Assume (i), and put $H=\mathcal H_r(T)$, $K=\mathcal K_r(T)$,
and $L=H+K$. By Proposition \ref{dichotomy}, $H\cap K=\{0\}$,
and $L$ is closed by Theorem \ref{mixed}(ii).
Since $H\subseteq L$, we have $\mathcal H_r(T_L)=H$.
The tails of any determining backward chain in $K$ also
lie in $K\subseteq L$; hence $\mathcal K_r(T_L)=K$.
Corollary \ref{canonical}, applied to $T_L$, gives
\[
L=H\oplus K\quad\text{topologically},\qquad
\sigma(T_L)\cap\mathbb T_r=\varnothing.
\]
Let $\pi:X\to X/L$ be the quotient map and $\widetilde T$
the operator induced by $T$. By Lemma \ref{saturation},
$\sigma(\widetilde T)\subseteq\mathbb T_r$. As $X/L$ is
finite-dimensional, a finite-dimensional complement of $L$
gives an upper triangular representation of $T$, and hence
\[
\sigma(T)=\sigma(T_L)\cup\sigma(\widetilde T),\qquad
\sigma(T)\cap\mathbb T_r=\sigma(\widetilde T).
\]
Let $P$ be the Riesz projection of $T$ for $\sigma(T)\cap\mathbb T_r$, taking $P=0$ if this set is empty.
The Riesz projection of $T_L$ for the same spectral set is zero,
whereas that of $\widetilde T$ is the identity. Since
$\pi T=\widetilde T\pi$, we obtain
\[
P|_L=0,\qquad \pi P=\pi.
\]
It follows that $\mathcal N(P)=L$ and
$\operatorname{rank}P=\codim L$. Thus, with
$F_r=\mathcal R(P)$,
\[
X=H\oplus F_r\oplus K
\]
topologically. The spectrum of $T_{F_r}$ is
$\sigma(T)\cap\mathbb T_r$, and $F_r$ is finite-dimensional.
This proves (ii), (iii), and \eqref{defectrank}.

Suppose (ii), and let $P$ be the sum of the Riesz projections
corresponding to $\sigma(T)\cap\mathbb T_r$. Put
$F_r=\mathcal R(P)$ and $Y=\mathcal N(P)$.
Then $F_r$ is finite-dimensional,
$\sigma(T_{F_r})\subseteq\mathbb T_r$, and
$\mathbb T_r\subset\rho(T_Y)$. By Corollary \ref{canonical},
\[
Y=\mathcal H_r(T_Y)\oplus\mathcal K_r(T_Y).
\]
On $F_r$ both radial subspaces are zero. The direct-sum
identities therefore give
$\mathcal H_r(T)=\mathcal H_r(T_Y)$ and
$\mathcal K_r(T)=\mathcal K_r(T_Y)$. Hence (iii) holds.

The implication (iii)$\Rightarrow$(i) is immediate. Under (iii), put $L=\mathcal H_r(T)\oplus\mathcal K_r(T)$. The operator on $X/L$ is similar to $T_{F_r}$, of spectral radius at most $r$. Projecting independent preimages to $X/L$ and applying Remark~\ref{indradremark} gives $\mathcal I_r(T)\subseteq L$. By Lemma~\ref{indfinitecodim}, $\mathcal I_r(T)=\mathcal I_r(T_L)$. Now $\mathcal H_r(T_L)=\mathcal H_r(T)$ and $\mathcal K_r(T_L)=\mathcal K_r(T)$, since every tail of a determining chain belongs to $\mathcal K_r(T)$. Proposition~\ref{indordered} on $L$ gives $\mathcal I_r(T)=\mathcal K_r(T)$ and hence (iv).
Finally, (iv)$\Rightarrow$(i) follows from
Lemma \ref{indfiniterigidity}.
\end{proof}

\section{Cantor--Bendixson derivatives and \texorpdfstring{$g_{\sigma_*}^{\alpha}$}{g-sigma-alpha}-invertibility}
Put $\sigma_*\in\{\sigma,\sigma_b,\sigma_d\}$. For a compact set $E\subset\C$, define
\[
\operatorname{acc}^0 E=E,\qquad
\operatorname{acc}^{\beta+1}E=\operatorname{acc}(\operatorname{acc}^{\beta}E),\qquad
\operatorname{acc}^{\lambda}E=\bigcap_{\beta<\lambda}\operatorname{acc}^{\beta}E
\]
when $\lambda$ is a nonzero limit ordinal. Let $\omega_1$ be the first uncountable ordinal.
The Cantor--Bendixson theorem identifies $\operatorname{acc}^{\omega_1}E$ with the perfect kernel of $E$, namely its set of condensation points (points whose every neighborhood contains uncountably many points of $E$). In particular, $\operatorname{acc}^{\omega_1}E=\operatorname{acc}^{\omega_1+1}E$, and this set is empty if and only if $E$ is at most countable.
If two closed sets coincide on an open set, then their Cantor--Bendixson
derivatives coincide there at every ordinal stage.
For every bounded operator $S$ and $\sigma_*\in\{\sigma_b,\sigma_d\}$,
$\sigma(S)\setminus\sigma_*(S)\subseteq\operatorname{iso}\sigma(S)$.
Thus $\sigma(S)$ is at most countable if and only if $\sigma_*(S)$ is at most countable.

For a successor ordinal $\alpha=\beta+1$, an operator $T$ is called $g_{\sigma_*}^{\alpha}$-invertible if there is a topological $T$-invariant decomposition $X=M\oplus N$ such that $T_M$ is invertible and
\[
\operatorname{acc}^{\beta}\sigma_*(T_N)\subset\{0\}.
\]
For a nonzero limit ordinal $\alpha$, this means $g_{\sigma_*}^{\gamma}$-invertibility for some successor ordinal $\gamma<\alpha$; for $\alpha=0$, it means $0\notin\sigma_*(T)$, as in \cite[Definition 2.4]{AIO}.
By \cite[Theorem 2.6]{AIO}, $g_{\sigma_*}^{\alpha}$-invertibility is equivalent to $0\notin\operatorname{acc}^{\alpha}\sigma_*(T)$. We characterize this condition for successor ordinals by canonical radial subspaces, allowing their intersection and the codimension of their sum to be finite.
\begin{theorem}\label{radialtransfinite}
Let $T\in\B(X)$, let $\beta$ be an ordinal, and let $\sigma_*\in\{\sigma,\sigma_b,\sigma_d\}$. The following assertions are equivalent:
\begin{enumerate}
\item $0\notin\operatorname{acc}^{\beta+1}\sigma_*(T)$;
\item there is $r>0$ such that $X=\mathcal{H}_r(T)\oplus \mathcal{K}_r(T)$ algebraically and
\begin{equation*}
\operatorname{acc}^{\beta}\sigma_*(T_{\mathcal{H}_r(T)})\subset\{0\};
\end{equation*}
\item there is $r>0$ such that $X=\mathcal{H}_r(T)\oplus \mathcal{K}_r(T)$ algebraically and
\begin{equation*}
\operatorname{acc}^{\beta}\sigma_*(\widetilde T_r)\subset\{0\},
\end{equation*}
where $\widetilde T_r:X/\mathcal{K}_r(T)\to X/\mathcal{K}_r(T)$ is defined by
\[
\widetilde T_r(x+\mathcal{K}_r(T))=Tx+\mathcal{K}_r(T).
\]
\item there is $t>0$ such that
\[
\dim(\mathcal{H}_t(T)\cap \mathcal{K}_t(T))<\infty,\qquad
\codim(\mathcal{H}_t(T)+\mathcal{K}_t(T))<\infty,
\]
and $\operatorname{acc}^{\beta}\sigma_*(T_{\mathcal{H}_t(T)})\subset\{0\}$;
\item there are $r,s>0$ such that
\[
\dim(\mathcal{H}_r(T)\cap \mathcal{K}_s(T))<\infty,\qquad
\codim(\mathcal{H}_r(T)+\mathcal{K}_s(T))<\infty,
\]
and $\operatorname{acc}^{\beta}\sigma_*(T_{\mathcal{H}_r(T)})\subset\{0\}$.
\end{enumerate}
In (ii) and (iii) the decomposition is topological, and $\widetilde T_r$ is similar to $T_{\mathcal{H}_r(T)}$. In (iv) and (v) the two radial subspaces are closed.
\end{theorem}
\begin{proof}
Suppose (i). There is $\varepsilon>0$ such that
\[
\operatorname{acc}^{\beta}\sigma_*(T)\cap\{\lambda:0<|\lambda|<\varepsilon\}=\varnothing.
\]
The set $\sigma_*(T)\cap\{\lambda:0<|\lambda|<\varepsilon\}$ is at most countable. Indeed, otherwise some compact annulus contained in this punctured disc would meet $\sigma_*(T)$ in an uncountable compact set. Such a set contains a nonempty perfect subset, which is contained in every Cantor--Bendixson derivative of $\sigma_*(T)$, a contradiction. The same is true for
$\sigma(T)\cap\{\lambda:0<|\lambda|<\varepsilon\}$. Choose $0<r<\varepsilon$ with $\mathbb T_r\subset\rho(T)$. By Corollary \ref{canonical},
\[
X=\mathcal{H}_r(T)\oplus \mathcal{K}_r(T).
\]
Put $A=T_{\mathcal{H}_r(T)}$ and $B=T_{\mathcal{K}_r(T)}$. Then $\sigma(B)\subset\{\lambda:|\lambda|>r\}$ and $\sigma(A)\subset\{\lambda:|\lambda|<r\}$. Since $\sigma_*(T)=\sigma_*(A)\cup\sigma_*(B)$, localization of Cantor--Bendixson derivatives in the open disc gives
\[
\operatorname{acc}^{\beta}\sigma_*(A)
=\operatorname{acc}^{\beta}\sigma_*(T)\cap\{\lambda:|\lambda|<r\}
\subset\{0\}.
\]
This proves (ii).

Conversely, suppose (ii). Theorem \ref{mixed} gives the topological decomposition $T=A\oplus B$, with $B$ invertible. Hence
\[
\operatorname{acc}^{\beta}\sigma_*(T)
=\operatorname{acc}^{\beta}\sigma_*(A)\cup
\operatorname{acc}^{\beta}\sigma_*(B).
\]
The first set is contained in $\{0\}$, and the second is separated from $0$. Thus $0\notin\operatorname{acc}^{\beta+1}\sigma_*(T)$.

Finally, the map
\[
J:\mathcal{H}_r(T)\longrightarrow X/\mathcal{K}_r(T),\qquad Jx=x+\mathcal{K}_r(T),
\]
is an isomorphism, and $J T_{\mathcal{H}_r(T)}=\widetilde T_rJ$. Therefore (ii) and (iii) are equivalent.

Clearly, (ii) implies (iv), and (iv) implies (v) by taking $r=s=t$. Suppose (v), and put $A=T_{\mathcal{H}_r(T)}$. By Lemma \ref{indfiniterigidity}, $\sigma_*(T)$ and $\sigma_*(A)$ agree on a punctured neighborhood of zero. Their Cantor--Bendixson derivatives therefore agree on the same punctured neighborhood. Since $\operatorname{acc}^{\beta}\sigma_*(A)\subset\{0\}$, we obtain $0\notin\operatorname{acc}^{\beta+1}\sigma_*(T)$.
\end{proof}

\smallskip
\noindent Theorem \ref{radialtransfinite}, applied to $T-\mu I$, characterizes every point $\mu\notin\operatorname{acc}^{\beta+1}\sigma_*(T)$ by the radial subspaces of $T-\mu I$. Moreover, for a nonzero limit ordinal $\alpha$, $g_{\sigma_*}^{\alpha}$-invertibility is characterized by the existence of a successor ordinal $\gamma=\beta+1<\alpha$ satisfying any of (ii)--(v) of Theorem \ref{radialtransfinite}.

\medskip
\noindent A bounded operator is called Riesz if its nonzero translates are Fredholm, and meromorphic if every nonzero spectral point is a pole of its resolvent. An operator $R$ is called zeroloid if $\operatorname{acc}\sigma(R)\subset\{0\}$.

\begin{corollary}\label{classical}
For $r>0$ such that $X=\mathcal H_r(T)\oplus\mathcal K_r(T)$ algebraically, put $A_r=T_{\mathcal H_r(T)}$. Then
\begin{enumerate}
\item $T$ is generalized Drazin invertible if and only if $A_r$ is quasinilpotent for some such $r$.
\item $T$ is generalized Drazin--Riesz invertible if and only if $A_r$ is Riesz for some such $r$.
\item $T$ is generalized Drazin--meromorphic invertible if and only if $A_r$ is meromorphic for some such $r$.
\item $T$ is $g_z$-invertible if and only if $A_r$ is zeroloid for some such $r$.
\end{enumerate}
The restriction may be replaced by the quotient $\widetilde T_r$ of Theorem~\ref{radialtransfinite}. Equivalently, in each case one may require, for some $r,s>0$,
\[
\dim(\mathcal{H}_r(T)\cap \mathcal{K}_s(T))<\infty,\qquad
\codim(\mathcal{H}_r(T)+\mathcal{K}_s(T))<\infty,
\]
with the same condition on $T_{\mathcal{H}_r(T)}$ as above. Here one may also impose $r=s$.
\end{corollary}
\begin{proof}
By the decomposition characterizations in \cite[Theorem 7.1]{Koliha}, \cite[Theorem 2.3]{ZC}, \cite{ZD} and \cite[Theorem 4.17]{GZ}, the four cases correspond to $(\beta,\sigma_*)=(0,\sigma),(0,\sigma_b),(0,\sigma_d)$ and $(1,\sigma)$ in \cite[Definition 2.4]{AIO}. Apply \cite[Theorem 2.6]{AIO} and Theorem~\ref{radialtransfinite} in each case.
\end{proof}
An operator is almost invertible if it admits a reducing decomposition into an invertible operator and an operator with at most countable spectrum \cite{AIO}.
\begin{corollary}\label{almostradial}
Let $T\in\B(X)$. The following assertions are equivalent:
\begin{enumerate}
\item $T$ is almost invertible;
\item there is $t>0$ with $X=\mathcal H_t(T)\oplus\mathcal K_t(T)$ algebraically and $\sigma(T_{\mathcal H_t(T)})$ at most countable;
\item there are $r,s>0$ with $\dim(\mathcal H_r(T)\cap\mathcal K_s(T))<\infty$, $\codim(\mathcal H_r(T)+\mathcal K_s(T))<\infty$ and $\sigma(T_{\mathcal H_r(T)})$ at most countable.
\end{enumerate}
In (ii) one may use $\mathcal H_r(T)\oplus\mathcal K_s(T)$ at distinct radii, in either order, and replace the restriction by the operator induced on $X/\mathcal K_s(T)$ (with $s=t$ at equal radii). In (iii) one may prescribe $r=s$, $r<s$ or $s<r$. The radial subspaces in (ii) and (iii) are closed. In both statements $\sigma$ may be replaced by $\sigma_b$ or $\sigma_d$.
\end{corollary}
\begin{proof}
By the Cantor--Bendixson theorem and \cite[Theorem 1.1]{AIO},
$T$ is almost invertible if and only if
$0\notin\operatorname{acc}^{\omega_1+1}\sigma(T)$.
Theorem \ref{radialtransfinite}, with $(\beta,\sigma_*)=(\omega_1,\sigma)$,
gives the equivalence of (i)--(iii), since a perfect subset of $\{0\}$ is empty.

Suppose (ii) holds at a radius $t>0$. By Theorem \ref{mixed},
$\mathbb T_t\subset\rho(T)$. Choose $0<a<t<b$ with
$\sigma(T)\cap\{\lambda:a\leq|\lambda|\leq b\}=\varnothing$.
Corollary \ref{canonical} gives
\[
\mathcal H_a(T)=\mathcal H_t(T)=\mathcal H_b(T),\qquad
\mathcal K_a(T)=\mathcal K_t(T)=\mathcal K_b(T).
\]
Taking $(r,s)=(a,b)$ or $(b,a)$ gives the required direct sums, and hence (iii). In any exact decomposition the quotient is similar to the inner restriction by Theorem~\ref{mixed}. For the inner restriction $A$, both $\sigma(A)\setminus\sigma_b(A)$ and $\sigma(A)\setminus\sigma_d(A)$ consist of isolated spectral points and are therefore at most countable.
\end{proof}

\begin{proposition}\label{almostindcores}
If $0\notin\operatorname{acc}^{\omega_1}\sigma(T)$, then
\[
\mathcal I_r(T)=\mathcal K_r(T)\quad(r>0),
\qquad \mathcal I(T)=\mathcal K(T).
\]
\end{proposition}
\begin{proof}
By Corollary \ref{almostradial}, write $T=A\oplus B$ topologically,
where $A$ is invertible and $\sigma(B)$ is at most countable.
For $A$ the equality follows from the injectivity of $A$.
Fix $r>0$ and $x\in\mathcal I_r(B)$. Choose
$0<c<r^{-1}$ and independent preimages $B^nx_n=x$
such that $\limsup_n\|x_n\|^{1/n}<c$.
Since $\sigma(B)$ is countable, some $t\in(r,c^{-1})$
satisfies $\mathbb T_t\subset\rho(B)$. Let $P$ be the inner
Riesz projection at radius $t$. Then
\[
Px=(B_{\mathcal R(P)})^nPx_n,\qquad
\rad(B_{\mathcal R(P)})<t<c^{-1}.
\]
Thus $Px\in\mathcal I_{c^{-1}}(B_{\mathcal R(P)})=\{0\}$
by Remark \ref{indradremark}. The restriction of $B$ to $\mathcal N(P)$ is
invertible with inverse of spectral radius less than $t^{-1}$.
Thus $x\in\mathcal K_r(B)$. The reverse inclusion is immediate,
and the direct-sum identities complete the proof. Taking unions
over $r>0$ gives the last assertion.
\end{proof}

For a bounded operator $S$ with at most countable spectrum,
\[
\dim\mathcal K(S)<\infty \;\Longrightarrow\; \#\sigma(S)<\infty,
\qquad
\mathcal K(S)=\{0\} \;\Longrightarrow\; \sigma(S)\subseteq\{0\}.
\]
Indeed, the isolated points are dense in a countable compact spectrum. For each nonzero isolated point $\lambda$, the corresponding Riesz subspace is nonzero and lies in $\mathcal K(S)$ because the restriction of $S$ to it is invertible. These subspaces are linearly independent.

\begin{corollary}\label{zerofinitedefects}
For $T\in\B(X)$ the following are equivalent:
\begin{enumerate}
\item $T$ is generalized Drazin invertible;
\item $\mathcal H_0(T)$ and $\mathcal K(T)$ have finite-dimensional intersection and their sum has finite codimension;
\item for some $s>0$ the same two defects are finite for $\mathcal H_0(T)$ and $\mathcal K_s(T)$;
\item for some $r>0$ the two defects are finite for $\mathcal H_r(T)$ and $\mathcal K(T)$, and $\sigma(T_{\mathcal H_r(T)})$ is at most countable;
\item for some $t>0$ the two defects are finite for $\mathcal H_t(T)$ and $\mathcal K_t(T)$, and $\sigma(T_{\mathcal H_t(T)})\subseteq\{0\}$.
\end{enumerate}
In (ii)--(v), one may equivalently require an algebraic direct sum equal to $X$ in place of the two finite defects. Such a sum is then topological. In the exact version of (iv) one may use the operator induced on $X/\mathcal K(T)$ instead of the restriction.
In (iv) the countability condition is unchanged for $\sigma_b$ and $\sigma_d$.
\end{corollary}
\begin{proof}
The equivalence of (i) and (ii) is \cite[Proposition 3.3]{AC}, and (i)$\Leftrightarrow$(v) follows from Theorem \ref{radialtransfinite} with $(\beta,\sigma_*)=(0,\sigma)$. If (i) holds, then $\mathcal{H}_t(T)=\mathcal{H}_0(T)$ and $\mathcal{K}_t(T)=\mathcal{K}(T)$ for all sufficiently small $t>0$, so (iii) and (iv) follow.

Suppose (iii). Proposition \ref{dichotomy} yields $\mathcal{H}_0(T)\cap \mathcal{K}_s(T)=\{0\}$. By Lemma \ref{indfiniterigidity}, the two summands are closed and $T_{\mathcal{K}_s(T)}$ is invertible. Since $T_{\mathcal{H}_0(T)}$ is quasinilpotent, \eqref{localmixed} gives $\sigma(T)\cap\{\lambda:0<|\lambda|<\varepsilon\}=\varnothing$ for some $\varepsilon>0$. Thus (i) holds.

Suppose (iv), and put $A=T_{\mathcal{H}_r(T)}$ and $E=\mathcal{H}_r(T)\cap \mathcal{K}(T)$. By Lemma \ref{indfiniterigidity}, $A$ is bounded on a Banach space, and $E$ is finite-dimensional. Every determining chain for an element of $\mathcal{K}(A)$ is also a determining chain for $T$, so
\[
\mathcal{K}(A)\subset E,\qquad \dim \mathcal{K}(A)<\infty.
\]
By the observation preceding the corollary,
$\sigma(A)$ is finite. The identity \eqref{localmixed} now gives $0\notin\operatorname{acc}\sigma(T)$, proving (i).

Under (i), for all sufficiently small $t>0$ the two summands are $\mathcal H_0(T)$ and $\mathcal K(T)$, whose direct sum is $X$. Conversely, a direct sum equal to $X$ has both defects zero. The quotient assertion follows by similarity, and the assertion about countability follows as in Corollary~\ref{almostradial}.
\end{proof}

\begin{corollary}\label{smallcircles}
The operator $T\in\B(X)$ is generalized Drazin--Riesz invertible if and only if there is $r>0$ such that, for every $0<s<r$,
\[
\dim(\mathcal{H}_s(T)\cap \mathcal{K}_s(T))<\infty,\qquad
\codim(\mathcal{H}_s(T)+\mathcal{K}_s(T))<\infty.
\]
\end{corollary}
\begin{proof}
By Theorem~\ref{finitecircle}, the two defects are finite at $s$ exactly when $\mathbb T_s\cap\sigma_b(T)=\varnothing$. Thus the stated condition is equivalent to $0\notin\operatorname{acc}\sigma_b(T)$, which characterizes generalized Drazin--Riesz invertibility by \cite[Theorem 2.6]{AIO}.
\end{proof}

\section{Algebraic invariant decompositions}
For an invariant subspace $M$, not necessarily closed, the spaces $\mathcal{H}_r(T_M)$ and $\mathcal{K}_r(T_M)$ are defined using the norm inherited from $X$.

\begin{proposition}\label{arbitrary}
Suppose $X=M\oplus N$ algebraically, where $M,N$ are $T$-invariant. Assume that, for some $r>0$,
\begin{enumerate}
\item $T_M$ is injective;
\item $\mathcal{K}_r(T_M)=M$;
\item $\mathcal{H}_r(T_N)=N$.
\end{enumerate}
Then $M=\mathcal{K}_r(T)$ and $N=\mathcal{H}_r(T)$. In particular, both summands are closed, $T_M$ is invertible, and
\[
\lim_{k\to\infty}\|(T_M)^{-k}\|^{1/k}<r^{-1},\qquad
\lim_{k\to\infty}\|T_N^k\|^{1/k}<r.
\]
\end{proposition}
\begin{proof}
Fix $|\lambda|=r$ and write $x=m+n$, where $m\in M$ and $n\in N$. Take a chain $(m_j)$ in $M$ with $m_0=m$, $Tm_{j+1}=m_j$ and $\limsup_j\|m_j\|^{1/j}<r^{-1}$. The chain bound and (iii) imply that the series
\[
u=\sum_{j\ge0}\lambda^j m_{j+1},\qquad
v=-\sum_{j\ge0}\lambda^{-j-1}T^jn
\]
converge in $X$. Telescoping gives $(T-\lambda I)u=m$ and $(T-\lambda I)v=n$, so $T-\lambda I$ is surjective.
If $Tx=\lambda x$, invariance and directness give $Tm=\lambda m$ and $Tn=\lambda n$. Injectivity of $T_M$ forces $m_j=\lambda^{-j}m$ for every $j$, so the strict bound on the chain implies $m=0$. Likewise (iii) gives $n=0$. Hence $\mathbb T_r\subset\rho(T)$.

By Corollary \ref{canonical}, $X=\mathcal{K}_r(T)\oplus \mathcal{H}_r(T)$ topologically. Conditions (ii) and (iii) give $M\subset \mathcal{K}_r(T)$ and $N\subset \mathcal{H}_r(T)$. Since $X=M\oplus N$, we obtain $M=\mathcal{K}_r(T)$ and $N=\mathcal{H}_r(T)$. The uniform bounds now follow from \eqref{mixedbounds}.
\end{proof}

\begin{corollary}\label{Mn}
Let $\beta$ be an ordinal and $\sigma_*\in\{\sigma,\sigma_b,\sigma_d\}$. Then $T$ is $g_{\sigma_*}^{\beta+1}$-invertible if and only if there exist $r>0$ and $T$-invariant subspaces $M,N$ such that
\begin{enumerate}
\item $X=M\oplus N$ algebraically;
\item $T_M$ is injective and $\mathcal{K}_r(T_M)=M$;
\item $\mathcal{H}_r(T_N)=N$;
\item $\operatorname{acc}^{\beta}\sigma_*(T_N)\subset\{0\}$.
\end{enumerate}
Under (i)--(iii), both $M$ and $N$ are closed, so $T_N$ in (iv) is a bounded operator on a Banach space. Moreover,
\[
M=\mathcal{K}_r(T),\qquad N=\mathcal{H}_r(T).
\]
Equivalently, (iv) may be replaced by $\operatorname{acc}^{\beta}\sigma_*(\widetilde T)\subset\{0\}$, where $\widetilde T$ is induced by $T$ on $X/M$.
\end{corollary}
\begin{proof}
By Proposition \ref{arbitrary}, (i)--(iii) give $M=\mathcal{K}_r(T)$ and $N=\mathcal{H}_r(T)$, with both subspaces closed. Theorem \ref{radialtransfinite} and \cite[Theorem 2.6]{AIO} now apply. Conversely, take $M=\mathcal{K}_r(T)$ and $N=\mathcal{H}_r(T)$ in that theorem; the spectral bounds \eqref{mixedbounds} give (ii) and (iii). Finally, $T_N$ is similar to $\widetilde T$.
\end{proof}


\begin{thebibliography}{99}
\bibitem{AC} Z. Aznay and F. Chabbabi, \emph{Algebraic Direct Sums and Generalized Drazin Invertibility}, arXiv:2610.18433 (2026).
\bibitem{ACPF} Z. Aznay and F. Chabbabi, \emph{On Mbekhta's conjecture and pseudo-Fredholm operators on Banach and Hilbert spaces}, arXiv:2610.23229v1 (2026).
\bibitem{ACZReg} Z. Aznay, F. Chabbabi and H. Zariouh, \emph{Pseudo-Fredholm operators form a regularity on Hilbert spaces}, arXiv:2611.04004v1 (2026).
\bibitem{AIO} Z. Aznay, A. Ouahab and H. Zariouh, Almost invertible operators, \emph{Results Math.} \textbf{79} (2024), Article 201, doi:10.1007/s00025-024-02226-z.
\bibitem{GCP} Z. Aznay, A. Ouahab and H. Zariouh, Generalization of some commutative perturbation results, \emph{Ann. Funct. Anal.} \textbf{14} (2023), Article 6, doi:10.1007/s43034-022-00231-3.
\bibitem{GZ} Z. Aznay, A. Ouahab and H. Zariouh, On the $g_z$-Kato decomposition and generalization of Koliha Drazin invertibility, \emph{Filomat} \textbf{37} (2023), 2087--2103.
\bibitem{Danes} J. Dane\v{s}, On local spectral radius, \emph{\v{C}asopis pro p\v{e}stov\'an\'i matematiky} \textbf{112} (1987), 177--187, doi:10.21136/CPM.1987.118306.
\bibitem{GMO} M. Gonz\'alez, M. Mbekhta and M. Oudghiri, On the isolated points of the surjective spectrum of a bounded operator, \emph{Proc. Amer. Math. Soc.} \textbf{136} (2008), 3521--3528.
\bibitem{Koliha} J. J. Koliha, A generalized Drazin inverse, \emph{Glasgow Math. J.} \textbf{38} (1996), 367--381.
\bibitem{KP} J. J. Koliha and P. W. Poon, Spectral sets II, \emph{Rend. Circ. Mat. Palermo (2)} \textbf{47} (1998), 293--310.
\bibitem{Mbekhta} M. Mbekhta, Sur la th\'eorie spectrale locale et limite des nilpotents, \emph{Proc. Amer. Math. Soc.} \textbf{110} (1990), 621--631.
\bibitem{Rudin} W. Rudin, \emph{Functional Analysis}, 2nd ed., McGraw--Hill, New York, 1991.
\bibitem{ZC} S. \v C. \v Zivkovi\'c-Zlatanovi\'c and M. D. Cvetkovi\'c, Generalized Kato--Riesz decomposition and generalized Drazin--Riesz invertible operators, \emph{Linear Multilinear Algebra} \textbf{65} (2017), 1171--1193.
\bibitem{ZD} S. \v C. \v Zivkovi\'c-Zlatanovi\'c and B. P. Duggal, Generalized Kato-meromorphic decomposition, generalized Drazin-meromorphic invertible operators and single-valued extension property, \emph{Banach J. Math. Anal.} \textbf{14} (2020), 894--914.
\end{thebibliography}
\end{document}